\documentclass[a4paper,fleqn]{cas-sc}

\usepackage[authoryear]{natbib}

\def\tsc#1{\csdef{#1}{\textsc{\lowercase{#1}}\xspace}}
\tsc{WGM}
\tsc{QE}

\usepackage{multirow}
\usepackage{colortbl}
\definecolor{lightgray}{gray}{0.85}
\usepackage{comment}
\usepackage{algorithm}
\usepackage{subfig}

\usepackage{setspace}
\usepackage{lineno}

\newcommand{\oldrev}{}

\newcommand{\change}{}

\begin{document}

\let\WriteBookmarks\relax
\def\floatpagepagefraction{1}
\def\textpagefraction{.001}



\title [mode = title]{{Signatures-and-sensitivity-based multi-criteria variational calibration for distributed hydrological modeling applied to Mediterranean floods}}  



%

\author[1]{Ngo Nghi Truyen Huynh}[orcid=0000-0001-5078-3865]\cormark[1]
\ead{ngo-nghi-truyen.huynh@inrae.fr}
\credit{research plan, conceptualization, software and numerical result, analysis, draft preparation, final redaction}

\author[1]{Pierre-André Garambois}
\ead{pierre-andre.Garambois@inrae.fr}
\credit{research methodology, conceptualization, analysis, draft preparation, final redaction, supervision, project administration, funding acquisition}

\author[1]{François Colleoni}
\ead{francois.colleoni@inrae.fr}
\credit{conceptualization, software and numerical results, analysis, final redaction}

\author[1]{Pierre Javelle}
\ead{pierre.javelle@inrae.fr}
\credit{research methodology, conceptualization, analysis, final redaction, funding acquisition}





\affiliation[1]{organization={INRAE, Aix Marseille Université, RECOVER},
            addressline={3275 Route Cézanne}, 
            city={Aix-en-Provence},
            citysep={}, 
            postcode={13182}, 
            country={France}}



\begin{abstract}
Classical calibration methods in hydrology typically rely on a single cost function computed on long-term streamflow series. Even when hydrological models achieve acceptable scores in NSE and KGE, imbalances can still arise between overall model performance and its ability to simulate flood events, particularly flash floods. \change{Multi-scale signatures, which refer to hydrological signatures computed at different temporal and/or spatial scales, and distributed flood modeling, which accounts for spatial variability in input variables and model parameters, are important concepts in hydrological modeling. In this study, the potential of using multi-scale signatures is explored to enhance multi-criteria calibration methods for spatially distributed flood modeling, which remains considerable challenges}. 
We present a novel signatures and sensitivity-based calibration \change{approach implemented into a variational data assimilation algorithm capable to deal with high dimensional spatially distributed hydrological optimization problems. It is tested on 141 flash flood prone }catchments mostly located in the French Mediterranean region. Our approach involves computing several signatures, including flood event signatures, using an automated flood segmentation algorithm. We select suitable signatures for constraining the model based on their global sensitivity with the input parameters \change{through global signature-based sensitivity analysis (GSSA)}. We then perform two multi-criteria calibration strategies using the selected signatures, including a single-objective optimization approach, which transforms the multi-criteria problem into a single-objective function, and a multi-objective optimization approach, which uses a simple additive weighting method to select an optimal solution from the Pareto set. Our results show significant improvements in both calibration and temporal validation metrics, especially for flood signatures, demonstrating the robustness and delicacy of our signatures-based calibration framework for enhancing flash flood forecasting systems.
\end{abstract}




\maketitle






\section{Introduction}

Numerical hydrological models are used extensively to simulate catchments responses to atmospheric signals and are a key component of floods forecasting systems where accuracy in terms of peak location, amplitude and timing is crucial. As a matter of facts,  hydrological models, whatever their complexity and spatialization, consist in more or less empirical representations of flows through watersheds compartments and contain parameters that cannot be inferred directly from the available observations but can only be meaningfully estimated through a calibration procedure (e.g. \citet{gupta2006model, vrugt2008treatment}). Such procedures aim to improve the model capability in reproducing the available observations of hydrological responses dynamics by optimizing model parameters.

Nevertheless, the whole construction process of a hydrological model is faced with the issue of equifinality: distinct model structures and/or parameter sets can lead to similar (in a sense to be defined) simulations. The equifinality concept has been popularized in hydrology by \citet{beven1993prophecy} while the issues of uncertainty in determining environmental model structures and estimating their parameters were known (e.g. \citet{beck1987water, yeh1986review}). For a given hydrological model structure, the calibration of its parameters is in general an ill-posed inverse problem with non unique solutions. As a consequence, the definition of an optimization algorithm and of a calibration metric is an essential modeling decision. Indeed, it determines how hydrological information is seen and learnt in the calibration process and it can substantially affect the quality and consistency of model simulations. 

\change{The present contribution studies the use of hydrological signatures for model calibration, especially for  optimization of spatial fields of distributed models in view to enhance flood modeling capabilities. Indeed, hydrological signatures are "quantitative metrics that describe statistical or dynamical properties of hydrologic data series, primarily streamflow" allowing to extract
"meaningful information about watershed processes" \citep{mcmillan2021review}. 
In this context, the following crucial ingredients are introduced along with relevant literature:}
\change{\begin{enumerate}[(i)]
    \item cost function definition particularly for emphasizing information from flood signatures;
    \item time-varying signatures computation which requires streamflow signal segmentation and analysis, as well as signatures selection based on their global sensitivity to model parameters which provides valuable guidance for the choice of calibration metrics;
    \item signature-based calibration approaches linked to multi-criteria calibration. 
\end{enumerate}}
\noindent \change{All these ingredients are integrated in our proposed framework to bring time-varying signatures into global sensitivity analysis and multi-criteria calibration.}
\subsection{\change{Calibration metrics}}
\change{To begin, the objective function definition is essential for an optimization procedure in the sense it is the metric that determines how is measured the misfit between model outputs and observations of the quantities of interest}. In hydrology, most calibration approaches attempt to optimize input parameters of a model such that they result in a minimal misfit between simulated and observed discharge. \change{However, since} no single metric can exhaustively represent this misfit, the calibration of a hydrological model is "inherently multi-objective" as remarked by \citet{gupta1998toward}. Several performance metrics have been proposed over the past decades in the literature for hydrological modeling. The classical quadratic Nash-Sutcliffe efficiency (NSE) \citet{nash1970river} (cf. Eq. \ref{eq:appendix_NSE} in Appendix \ref{append:metrics}) has been used for long time. The Kling--Gupta (KGE) (cf. Eq. \ref{eq:appendix_KGE} in Appendix \ref{append:metrics}) proposed in \citet{gupta2009decomposition} and based on a decomposition of the NSE has also become widely used. 
Other metrics, in form of signature measures (see review in \citet{mcmillan2021review}), have been proposed in the literature for model evaluation (e.g. \citet{Yilmaz2007WR006716}) and used in model optimization (e.g. \change{\citet{roux2011physically,Shafii2014WR016520,mostafaie2018comparing, Kavetski_2017WR, sahraei2020signature,wu2021improved}} and references therein). Hydrological signatures can be used to derive application-specific metrics such as for high flows in \citet{mizukami2019choice,roux2011physically}.
\subsection{\change{Hydrological signatures and global sensitivity analysis}}
\change{Moreover,} hydrological signatures are a useful tool to effectively evaluate models and diagnose the role of their components in explaining the discrepancy between the simulated and observed behavior \citep{gupta2009decomposition}, especially when \change{signatures are combined with global sensitivity analysis \citep{horner_PhD_2020}. Indeed sensitivity analysis (SA) examines how the variation of a model output, and consequently of a simulated signature, can be apportioned to a variation in its inputs \cite{saltelli2002making}. SA enables to establish which parameters in a model most importantly affect the magnitude, variability and dynamics of model response \citep{Razavi_Gupta_SA_2014WR}, to identify signatures-parameters links \citep{horner_PhD_2020}. In contrast to local SA, which focuses on a specific point in the model parameter space, global sensitivity analysis (GSA) considers the whole variation range of the inputs \citep{Saltelli_2008}. GSA have  been developed in a statistical framework (see review in \cite{iooss2015review}) and extensively applied in hydrological modeling (see review in \cite{SONG2015}; see also \cite{gupta2018revisiting,razavi2019multi} and references therein). Efficient GSA methods for estimating Sobol' sensitivity indices \citep{sobol1990sensitivity,Saltelli_2008} and libraries are now available (see for example a recent benchmark in \cite{puy2022comprehensive}). Although using GSA to perform signature-based sensitivity analysis is an interesting topic, it remains poorly studied. 
Nevertheless, would it be for model diagnostic, sensitivity analysis or multi-criteria calibration, the computation of hydrological signatures at varying time scales is faced with the difficulty of consistent segmentation of flood events.
This issue has been highlighted by \citet{tarasova2018exploring} and will be detailed later. In the present work, an original segmentation algorithm is proposed on top of a global signature-based sensitivity analysis (GSSA).}
\subsection{\change{Multi-criteria calibration approaches}}
\change{Last but not least, the definition of the calibration algorithm itself is crucial and has been extensively studied in hydrological modeling.} Hydrological calibration problems that incorporate multiple metrics, including multi-scale signatures, \change{which we define as signatures at different temporal and/or spatial scales,} can be considered as multi-criteria optimization problems \citep{gupta1998toward}. Generally, three categories of methods are employed for solving multi-criteria optimization problems in various domains: (i) transforming the multi-criteria problem into a single-objective optimization problem \citep{ross2015multiobjective, el2014optimal, veluscek2015composite}; (ii) obtaining a non-inferior solution set (Pareto front) by solving the multi-objective optimization problem \citep{khorram2014numerical, tavakkoli2011new, torres2011multi}; (iii) selecting a unique solution after obtaining the Pareto optimal solution set by adding constraints based on specific preferences \citep{chibeles2016multi, wu2015pareto}.
The state-of-the-art in multi-criteria optimization in hydrology is commonly accomplished through the first two approaches mentioned earlier. For instance, a simple approach on the choice of \change{single} calibration metric for flood modeling, including NSE, \change{or} weighted KGEs, \change{or} annual peak flow bias, has been proposed for daily mHm and VIC models on 492 US catchments by \citet{mizukami2019choice}. For event-based flash flood modeling at high resolution, a metric that, \change{in addition to NSE, } accounts for the shape of flash flood hydrographs, particularly their timing and maximum peak flow, has been studied in \citet{roux2011physically}. \change{These two studies in a flood modeling context, highlight the interest of signature-based calibration and the need of further investigations. They only consider a limited set of flood signatures, that is chosen empirically and computed without an automated segmentation algorithm. 
Furthermore, existing methods for signature-based calibration are only capable to solve low dimensional optimization problems which remains a significant limitation, especially to tackle (regional) spatially distributed modeling and high dimensional optimization problems with multi-source data (e.g. hundreds or more tunable parameters) for which variational data assimilation (VDA) approach is well suited (cf. \cite{jay2020potential})}. \change{Despite the potential benefits of using multi-scale hydrological signatures in calibration,} generalizing these methods and integrating them into VDA algorithms remain significant challenges\change{. This is mainly} due to the complexity involved in consistently computing these signatures while ensuring differentiability \change{for high dimensional context}. 
\change{In addition to the aforementioned methods for integrating multi-criteria problem into a single calibration metric}, research on calibration with multi-objective functions to generate a set of non-dominated solutions has also been conducted, as seen in studies by \citet{yapo1998multi, guo2014multi, oliveira2021contribution, mostafaie2018comparing}. However, the selection of an optimal solution from the non-dominated set has not received significant attention. 
Our goal in this work is to comprehensively investigate \change{signature-based calibration, guided by signature-parameters links estimated through GSSA, with} all feasible multi-criteria optimization methods\change{, including an adjoint based VDA, to improve generalizability.}

\change{To summarize, this research will study four aspects that have received relatively little attention in prior studies: 
\begin{enumerate}[(i)]
    \item the need for an automated segmentation method applicable to large contrasted catchment-floods samples and capable to capture hydrological information at the scale of flash flood events;
    \item a global analysis of simulated errors across various hydrological signatures and their sensitivity with the model parameters;
    \item the need for a more intelligent approach to select the Pareto optimal solution in the case of optimization with multi-objective functions;
    \item the computation of the cost function based on signatures within a VDA algorithm capable to deal with large  spatialized parameter vectors.
\end{enumerate}}
\noindent The proposed framework originally integrates automated segmentation of flood events and signatures computation within a VDA algorithm from \citet{jay2020potential}, enabling high dimensional spatially distributed calibration with multi-criteria metrics adapted to flood modeling. Classical global calibration algorithms have also been upgraded that way. These upgrades, including new cost functions and adjoint model update, have been implemented into the SMASH platform, which solvers are differentiable. 
Using the proposed algorithms, we investigate over a quite large dataset of Mediterranean flash floods the parametric sensitivity of a parsimonious distributed hydrological model for a large array of signatures from the literature, as well as the benefit of using a signature-based flood specific metric in calibration, and especially in performing variational spatially distributed optimization which has seldom been done to our best knowledge.

The remaining sections of this paper are organized as follows: section \ref{sec:Methodology} describes our methodology for computing various hydrological signatures and our multi-criteria calibration algorithms, along with an overview of the SMASH forward model. In section \ref{sec:Result}, we present and analyze our results on signatures and calibration, including a summary of the data and numerical experiments. Finally, in section \ref{sec:Conclusion}, we conclude our work and outline potential future directions.

\section{Methodology\label{sec:Methodology}}

\begin{figure}[ht!]
\begin{centering}
\includegraphics[scale=0.78]{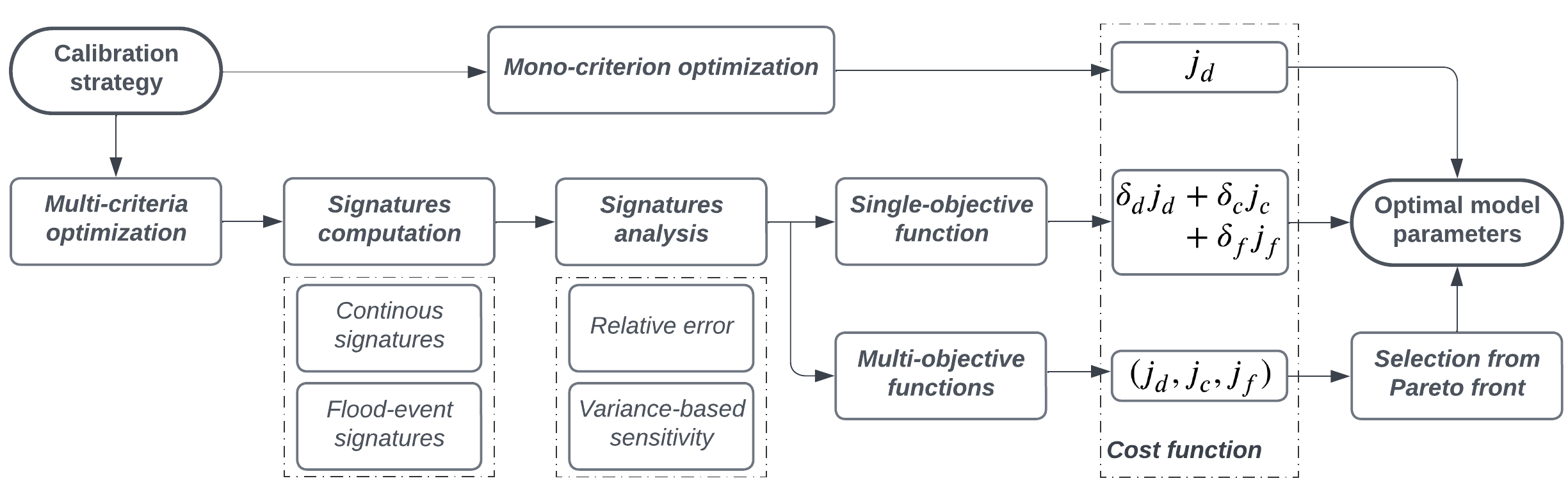}
\par\end{centering}
\centering{}\caption{\oldrev{Flowchart of the multi-criteria calibration process using hydrological signatures. The different cost functions are denoted by $j_d$, $j_c$ and $j_f$, while the corresponding optimal weights are denoted by $\delta_d$, $\delta_c$ and $\delta_f$. The notations used in the cost function will be explained in \ref{subsec:signatures-based-multi-criteria}.}}
\label{calib_flowchart}
\end{figure}

\oldrev{We propose a novel calibration strategy that leverages hydrological signatures and their sensitivity analysis in combination with the optimization algorithms discussed above. Our approach is illustrated in Fig. \ref{calib_flowchart} and addresses the challenges of model calibration in the presence of multiple objectives and complex hydrological processes.}

\oldrev{The computations of the signatures are first performed to quantify their sensitivities with the model parameters following \citet{horner_PhD_2020}. These computations involve performing both whole-period-based analysis to obtain continuous signatures and event-based analysis to capture the most significant events (flood event signatures). Through this analysis, we gain a more meaningful understanding of the parametric sensitivity, not just for discharge but also for other factors that need to be considered as part of our minimization criterion. Furthermore, we evaluate the sensitivity of signature error using variance-based sensitivity analysis (Sobol indices) to determine the most appropriate signatures for multi-criteria optimization. Based on these results, we conduct a multi-criteria optimization with single-objective or multi-objective functions, utilizing suitable hydrological signatures to improve the simulation performance.}

The numerical algorithms proposed here are implemented in Python, on top of SMASH Fortran platform that is interfaced in Python \citep{JayCoGa_2022_wrappingSMASH} making accessible its forward-inverse algorithms (forward hydrological models, \change{Step-By-Step} (SBS\change{, steepest descent algorithm summarized in \citet{edijatno91}}) and VDA \change{\citep{jay2020potential}} calibration algorithms) and internal variables.

\oldrev{The following subsections of this section detail the different elements of our methodology: \ref{subsec:forward_and_calibration_algos} defines the hydrological model structure, the objective function and the proposed calibration algorithms; \ref{subsec:signatures computation} explains which signatures are computed and how, including a description of the proposed hydrograph segmentation algorithm; \ref{subsec:signatures sensitivity} describes the method for computing global sensitivities of simulated hydrological signatures; \ref{subsec:signatures-based-multi-criteria} details the formulation of the multi-criteria cost functions including multi-scale signatures and the multi-objective optimization problems.}

\subsection{\oldrev{SMASH: An overview of the forward model and calibration algorithms}}
\label{subsec:forward_and_calibration_algos}

SMASH is a computational software framework dedicated to \textit{Spatially distributed Modelling and ASsimilation for Hydrology}. It aims to tackle flexible spatially distributed hydrological modeling, signatures and sensitivity analysis, as well as high dimensional inverse problems using multi-source observations. This model is designed to simulate discharge hydrographs and hydrological states at any spatial location within a basin and reproduce the hydrological response of contrasted catchments, especially \change{aimed} at floods and low-flows modeling, by taking advantage of spatially distributed meteorological forcings, physiographic data and hydrometric observations.

First, the forward spatially distributed hydrological modeling problem is formulated as follows. Let $\Omega\subset\mathbb{R}^{2}$ be a 2D spatial domain (catchment) and $t>0$ be the physical time. A regular lattice $\mathcal{R}_{\Omega}$ covers $\Omega$ and $D(x)$ is the drainage plan obtained from terrain elevation processing. 
The number of active cells within a catchment $\Omega$ is denoted $N_{x}$. 
Then the hydrological model is a dynamic operator $\mathcal{M}$ mapping observed input fields of rainfall and evapotranspiration  $\boldsymbol{P}\left(x,t'\right)$,  $\,\boldsymbol{E}\left(x,t'\right)$,  $\forall\left(x,t'\right)\in\Omega\times[0,t]$ onto discharge field $Q(x,t)$ such that:
\begin{equation}\label{eq:forward model}
    Q\left(x,t\right) = \mathcal{M}\left[\boldsymbol{P}\left(x,t'\right),\boldsymbol{E}\left(x,t'\right),\boldsymbol{h}\left(x,0\right),\boldsymbol{\theta}\left(x\right),t\right], 
    \forall x\in\Omega, t'\in\left[0,t\right]
\end{equation}
with $\boldsymbol{h}\left(x,t\right)$ the $N_{s}$-dimensional vector of model states 2D fields and $\boldsymbol{\theta}$ the $N_{p}$-dimensional vector of model parameters 2D fields. In the following, $\boldsymbol{\theta}$ is also called control vector in optimization context.

\begin{figure}[ht!]
     \centering
     \includegraphics[scale=0.25]{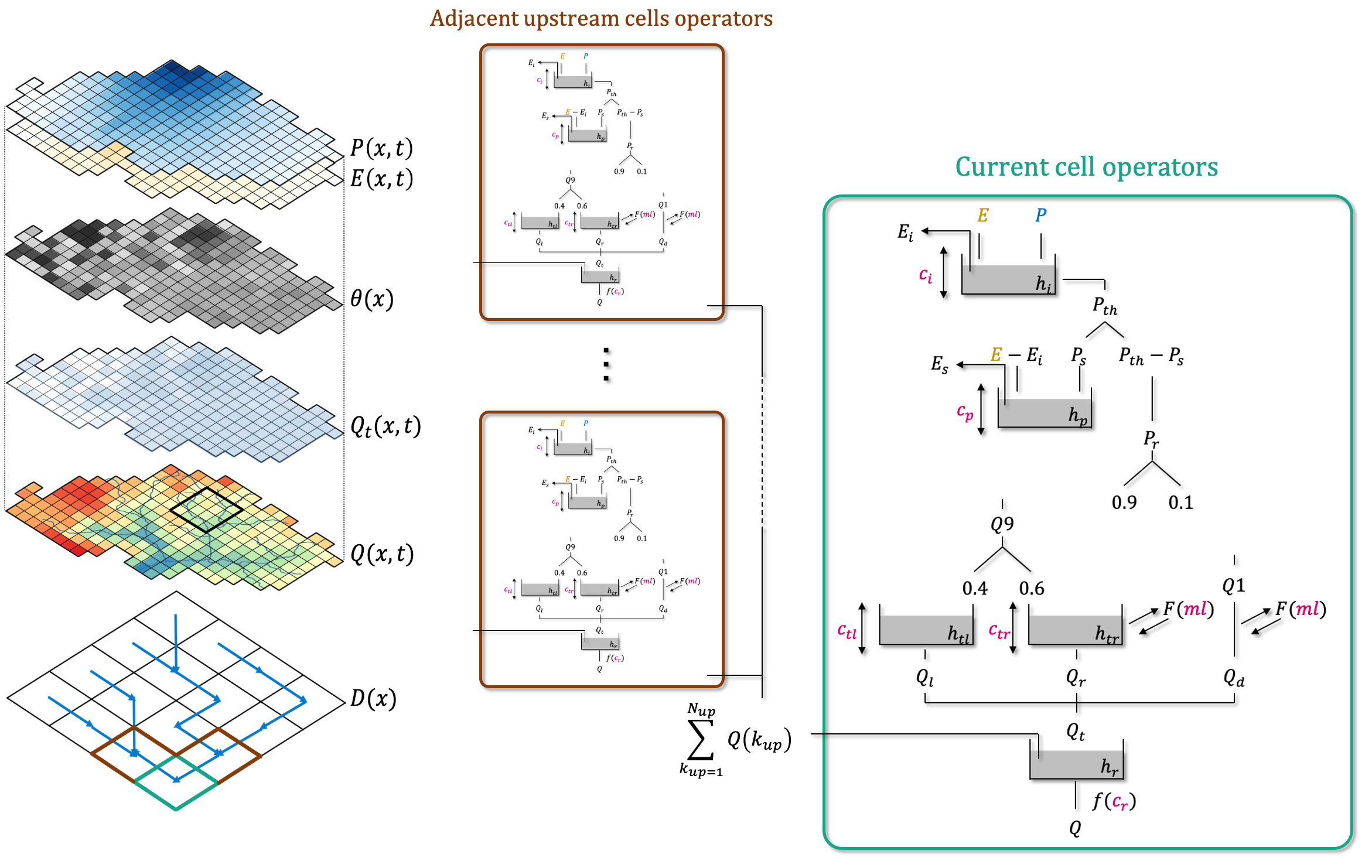}
     \caption{Distributed hydrological modeling with SMASH platform. \oldrev{Model fields from top to bottom: meteorological inputs, parameters, internal and output flux maps (left). Pixel scale and pixel-to-pixel flow operators of SMASH model structure S6 studied (right)}. 
     }
     \label{fig:SMASH_grid_model}
 \end{figure}

Then, the forward hydrological model structure with the parsimonious 6 parameters from \citet{egusphere-2022-506} 
(Fig. \ref{fig:SMASH_grid_model}) is defined as follows. 
For a given cell $i$ of coordinates $x\in\Omega$, in the proposed model S6, four reservoirs $\mathcal{I}$, $\mathcal{P}$, $\mathcal{T}_r$ and $\mathcal{T}_l$ of respective capacity $c_{i}$, $c_{p}$, $c_{tr}$ and  $c_{tl}$, are considered for simulating, respectively, the interception, the production of runoff and its transfer within a cell. 
Their state vector is denoted
$\boldsymbol{h}(x,t)\equiv\left(\boldsymbol{h}_{i}(x,t),\boldsymbol{h}_{p}(x,t),\boldsymbol{h}_{tr}(x,t),\boldsymbol{h}_{r}(x,t),\boldsymbol{h}_{tl}(x,t)\right)^{T}$,
and the parameter vector of SMASH model structure S6 is
$\boldsymbol{\theta}(x)\equiv\left(c_{i}(x),c_{p}(x),c_{tr}(x),c_{r}(x),ml(x),c_{tl}(x)\right)^{T}$.
Hence the size of state vector is $N_{s}\times N_{x}$ and the size of parameter vector that is optimized in the following is $N_{p}\times N_{x}$. Considering tens of cells or more over a simulated catchments domain $\Omega$, the calibration of $\boldsymbol{\theta}$ is a high dimensional inverse problem. 
The numerical model operates at hourly time step $dt=1h$ and on a regular grid at $dx=1km$.

In order to calibrate the hydrological model based on the simulated and observed discharge at gauged cells $x_{k}\in\Omega$, $k\in 1,..,Ng$, denoted as $Q_k(t)$ and $Q^*_k(t)$, respectively, we define the objective convex function as shown in Eq. \ref{eq:cost function}.
\begin{equation}\label{eq:cost function}
    J(\boldsymbol{\theta}) = J_{obs} (\boldsymbol{\theta}) + \alpha J_{reg}(\boldsymbol{\theta})
\end{equation}
where the observation cost function 
$J_{obs}=\frac{1}{N_{g}}\sum_{k=1}^{N_{g}}\text{\ensuremath{j^*_{k}}}$
measuring the misfit, via several adapted metrics \oldrev{that can include signatures as} detailed later, between simulated and observed discharge. In this study, $N_{g}=1$, that is for single gauge optimization. 
Note that simulated discharge
$Q_k(t) = \mathcal{M}\left[\boldsymbol{P}\left(x,t'\right),\boldsymbol{E}\left(x,t'\right),\boldsymbol{h}\left(x,0\right),\boldsymbol{\theta}\left(x\right)\right],
\forall x\in\Omega_k, t'\in\left[0,t\right]$
with $\Omega_k\subset\Omega$ denoting the spatial domain including all upstream cells of a gauge at $x_k$, depends on the control vector $\boldsymbol{\theta}$ via the hydrological model $\mathcal{M}$ (Eq. \ref{eq:forward model}). The second term in Eq. (\ref{eq:cost function}) is weighted by $\alpha$ and set as a classical Thikhonov regularization
$J_{reg}=\left\Vert \boldsymbol{B}^{-1/2}\left(\boldsymbol{\theta}-\boldsymbol{\theta}^{*}\right)\right\Vert _{L^{2}}^{2}$
with $\boldsymbol{B}$ the background error covariance, and $\boldsymbol{\theta}^*$ the first guess/background on $\boldsymbol{\theta}$. We set $\alpha = 10^{-4}$ for the spatially distributed optimizations presented in this study, $\alpha = 0$ otherwise if $\boldsymbol{\theta} \equiv \overline{\boldsymbol{\theta}}$, and $\boldsymbol{B}$ is simply defined from $\boldsymbol{\sigma}_{\theta}$ the vector of mean deviations of $\boldsymbol{\theta}$, as done in \citet{jay2020potential}. 
\oldrev{The optimal estimate $\hat{\boldsymbol{\theta}}$ of the model parameter set can be obtained by minimizing the objective function $J$ in Eq. \ref{eq:forward model}, subject to an additional bound constrain on the model parameters, which can be expressed as Eq. \ref{inv_problem}.}
\begin{equation}
\hat{\boldsymbol{\theta}}=\underset{\boldsymbol{\theta}_{min}\leq\boldsymbol{\theta}\leq\boldsymbol{\theta}_{max}}{\arg\min}J\left(\boldsymbol{\theta}\right) \label{inv_problem}
\end{equation}
This inverse problem \ref{inv_problem} is tackled with different global optimization algorithms considering a spatially uniform control, that is low dimensional optimization problems. For instance, optimization algorithms such as: Step-By-Step (SBS), Nelder–Mead and Genetic Algorithms (GA) can be applied in this scenario. Next, a spatially distributed control vector is sought with a VDA algorithm \citep{jay2020potential} adapted to such high dimensional hydrological optimization problems.
Considering a spatially distributed control vector $\boldsymbol{\theta}(x)$, its optimization is performed with the L-BFGS-B algorithm (limited-memory Broyden–Fletcher–Goldfarb–Shanno bound-constrained \citep{Zhu1997}) adapted to high dimension. This algorithm requires the gradient of the cost function with respect to the sought parameters  $\nabla_{\boldsymbol{\theta}} J $, that is obtained by solving the adjoint model. \change{The} numerical adjoint model has been generated with the automatic differentiation engine TAPENADE \citep{tapenade} applied to the SMASH source code\change{, including the novelties added into the forward code,} and validated with standard gradient test.
The background value $\boldsymbol{\theta}^{*}$, used as a starting point for the optimization problem and in the regularization term, is set as in \citet{jay2020potential}, i.e. as $\bar{\boldsymbol{\theta}}$, a spatially uniform global optimum determined with a simple global-minimization algorithm from a uniform first guess $\bar{\boldsymbol{\theta}}^*$. Given mildly non linear hydrological models as those considered in this study, this calibration approach is pertinent and sensitivity to priors is limited as shown in \citet{jay2020potential}.

\subsection{Signatures computation}\label{subsec:signatures computation}

\begin{figure}[ht!]
\begin{centering}
\includegraphics[scale=0.6]{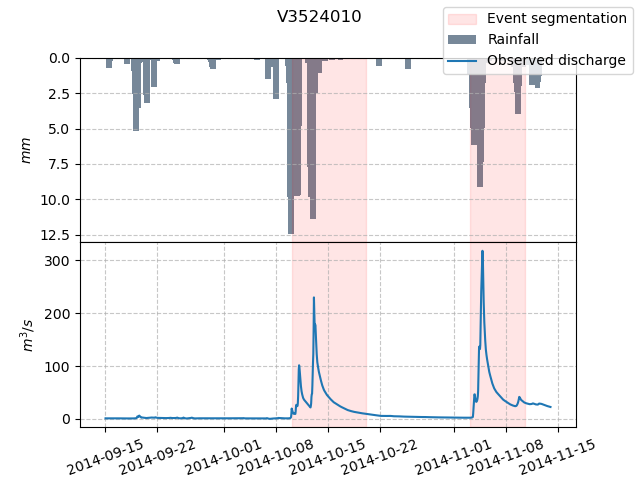}
\par\end{centering}
\centering{}\caption{\oldrev{Example of flood events detected from hydrograph using the segmentation algorithm.}}
\label{hydrograph_event}
\end{figure}

Several signatures describing and quantifying properties of discharge time series are introduced in view to analyze and calibrate hydrological models \oldrev{(an exhaustive list is given in Appendix \ref{appendix_sign})}. Signatures are denoted as $S_{i}, \ i\in{1..N_{crit}}$, with $N_{crit}$ being the number of different signature types considered. \oldrev{These signatures allow for the description of various aspects of the rainfall--runoff behavior, such as flow distribution (e.g. based on flow percentiles), flow dynamics \citep{lemesnil:tel-03578569}, flow separation \citep{nathan1990evaluation,lyne1979stochastic}, and flow timing, among others.} A so-called continuous signature is a signature that can be computed \oldrev{over the entire study period}.
Flood event signatures on the other hand focus on the behavior of the high-flows that are observed in flash flood events (Fig. \ref{hydrograph_event}).

\change{The computation  of multi-scale signatures from hydrological time series, typically at flood event scale, requires a segmentation approach.
Although the concept of flood event is widely used in hydrology, there is no clear consensus on approaches for flood detection from continuous streamflow time series, as pointed out in \citet{tarasova2018exploring}.
Several studies have suggested segmentation algorithms for detecting flood events (refer to the references in \citet{tarasova2018exploring}). For instance, \citet{li2022development, astagneau_2021} used simple segmentation methods respectively involving fixed time windows before and after rainfall events or discharge thresholds to detect events. 
Meanwhile, \citet{tarasova2018exploring} developed an algorithm incorporating, baseflow separation technique (see also \citet{pelletier2020hydrograph}), rainfall attribution methods and an iterative procedure to identify single-peak components of multiple-peak events. In this study, we propose an automated segmentation algorithm, consisting of, peak detection in discharge series, catchment rainfall time series analysis through a combination of rainfall gradients and rainfall energy criterion, which enables a robust determination of flood start time on contrasted catchment-floods, and a classical baseflow separation for determining the end of an event (see Appendix \ref{app:sec:seg-algo} for a detailed explanation of our algorithm).}

\subsection{Signatures sensitivity}\label{subsec:signatures sensitivity}

To perform a calibration process with hydrological signatures, it is important to investigate the sensitivity of simulated signatures with the model parameters, to guide the potential selection of the signatures which should be used to calibrate the model. The sensitivity analysis enables us to examine how the variation of a given output/signature can be apportioned to a variation in model inputs \citep{saltelli2002making}. If \change{certain} signatures are not sensitive with the model parameters, \oldrev{then it may not have any significant impact} to optimize an objective function based on these signatures.
\oldrev{In this context, we consider a hydrological model $\mathcal{M}$ with $m$ spatially uniform parameters $\bar{\boldsymbol{\theta}} \equiv (\boldsymbol{\theta}_{1},...,\boldsymbol{\theta}_{m})$. Then the simulated value of a signature $S_i$, calculated from the simulated discharges via a discharge-to-signature mapping $f_i$, is represented as $S_{i}^{s} \equiv f_i \circ \mathcal{M}\left(\boldsymbol{P},\boldsymbol{E},\boldsymbol{h}, \bar{\boldsymbol{\theta}}\right)$.} We are interested in Sobol indices called first-order \oldrev{and total-order. }
The first- (depending on $\boldsymbol{\theta}_{j}$) and \oldrev{total- (depending on $\boldsymbol{\theta}_{\sim j}$, i.e. all parameters except $\boldsymbol{\theta}_{j}$)} Sobol indices \oldrev{of the simulated signature} $S_{i}^s$ are respectively defined as follows:
\begin{equation*}
    s_{i}^{(1j)}=\frac{\mathbb{\mathbb{V}}[\mathbb{E}[S_{i}^s|\boldsymbol{\theta}_{j}]]}{\mathbb{\mathbb{V}}[S_{i}^s]}\oldrev{=\frac{V_j}{V} \;\;\;\;}
    \text{ and \;\;\;}
    \oldrev{s_{i}^{(Tj)}=\frac{\mathbb{\mathbb{E}}[\mathbb{V}[S_{i}^s|\boldsymbol{\theta}_{\sim j}]]}{\mathbb{\mathbb{V}}[S_{i}^s]}=1-\frac{\mathbb{\mathbb{V}}[\mathbb{E}[S_{i}^s|\boldsymbol{\theta}_{\sim j}]]}{\mathbb{\mathbb{V}}[S_{i}^s]}=1-\frac{V_{\sim j}}{V}}
\end{equation*}
\oldrev{where $V_j$ (respectively, $V_{\sim j}$) is the variance of the expectation of output signature $S_{i}^s$ conditioned  by the input parameter $\boldsymbol{\theta}_j$ (respectively, $\boldsymbol{\theta}_{\sim j}$, i.e. all sampled inputs except  $\boldsymbol{\theta}_j$).
To estimate these indices, \citet{azzini2021comparison} proposed a method based on the Saltelli generator \citep{saltelli2002making}, which is implemented in the \textit{SALib} Python library \citep{Iwanaga2022, Herman2017}. This method, that is shown to be relatively accurate in a recent benchmark \citep{puy2022comprehensive}, allows us to estimate the first-, second- and total-order variance-based sensitivity indices using Monte Carlo simulations. However, in our specific application with high dimensional parameter spaces, we have encountered significant challenges in estimating the second-order variance-based sensitivity indices due to their computationally intensive nature \citep{saltelli2002making, campolongo2011screening}. To achieve accurate results, a large number of Monte Carlo simulations are required, which can be time-consuming and computationally demanding. Therefore, for the purpose of this study, we focus on estimating the first- and total-order Sobol indices, which provide a sufficiently efficient means of capturing information about interaction effects while retaining an acceptable computational cost.}

\subsection{\oldrev{Multi-criteria calibration using hydrological signatures}}
\label{subsec:signatures-based-multi-criteria}
\oldrev{This section defines the calibration objective functions and how they account for the multi-scale signatures that are provided by the segmentation algorithm detailed previously.}

\oldrev{First, we define cost function parts corresponding respectively to classical metrics, continuous signatures and event based signatures.} Let us consider a classical objective function $j_{d}$, which is the dominant criterion (or the most constrained criterion) in case of multi-criteria optimization, an objective function $j_{c}$
combining continuous-signatures-based cost functions, and $j_{f}$
combining flood-event-signatures-based cost functions. Then, the
cost function to be minimized, denoted $J$, can be defined as Eq. \ref{eq:totalcost}.
\begin{equation}
J\equiv
    \begin{cases}
        \delta_{d}j_{d}+\delta_{c}j_{c}+\delta_{f}j_{f} \text{ for single-objective optimization,}\\
        (j_{d},j_{c},j_{f}) \text{ for multi-objective optimization}
    \end{cases}
\label{eq:totalcost}
\end{equation}
where $\delta_{d},\delta_{c},\delta_{f}$ are the corresponding optimization
weights in the first case. Keep in mind that we take into account the use of signatures
in both cases but the first case is a single-objective optimization
while the second is a multi-objective optimization.

\oldrev{Then we detail how each cost function part is computed from signatures.} For each signature $S_{i}$, denote by $S_{i}^{o}$ and $S_{i}^{s}$
 the observation and simulation respectively. The set of continuous
and flood event signatures denoted $N_{c}$ and $N_{f}$ respectively.
Then, the components $j_{d}$, $j_{c}$ and $j_{f}$ can be defined
as follows:
\begin{itemize}
\item $j_{d}\equiv1-NSE$ or $1-KGE_{\left(\alpha,\beta,\gamma\right)}$
with varying weights $\alpha,\beta,\gamma$ (see Appendix \ref{append:metrics}). This metric $j_d$ is considered
as a constraining objective function for selecting an optimal solution from non-inferior solutions in case of multi-objective optimization (see Appendix \ref{select_sol_append}).
\item $j_{c}\equiv$
    $\begin{cases}
        \sum_{S_{i}\in N_{c}}\sigma_{S_{i}}j_{c}^{S_{i}}\text{, for single-objective multi-criteria optimization;}\\ 
        (\{j_{c}^{S_{i}}\}_{S_{i}\in N_{c}})\text{, for multi-objective optimization}
    \end{cases}$
    
where $j_{c}^{S_{i}}\equiv\left|\frac{S_{i}^{s}}{S_{i}^{o}}-1\right|$
is the objective function based on continuous signature $S_{i}$ and
$\sigma_{S_{i}}$ is the corresponding optimization weight of $S_{i}$ in case of single-objective function.
\item $j_{f}\equiv$
    $\begin{cases}
        \sum_{S_{i}\in N_{f}}\sigma_{S_{i}}j_{f}^{S_{i}}\text{, for single-objective multi-criteria optimization;}\\
        (\{j_{f}^{S_{i}}\}_{S_{i}\in N_{f}})\text{, for multi-objective optimization.}
    \end{cases}$

In this case, and in the context of global optimization in time,
$j_{f}^{S_{i}}\equiv\frac{1}{N_E}\sum_{e=1}^{N_{E}}\left|\frac{S_{i,e}^{s}}{S_{i,e}^{o}}-1\right|$ defines the scalar objective function related to flood signature ${S_{i}\in N_{f}}$ over the $N_{E}$ events selected with the segmentation method described in Algorithm \ref{algo_seg}.
Otherwise, to perform a season-based optimization on flood event signatures, we can compute for the events occurring in the selected season. For example, for a Spring-based optimization: 
\begin{equation*}
    j_{f,spring}^{S_{i}}\equiv\frac{1}{\dim \mathcal{SP}}\sum_{e\in\mathcal{SP}}\left|\frac{S_{i,e}^{s}}{S_{i,e}^{o}}-1\right| 
    \text{ s.t. } 
    \forall e\in\mathcal{SP}\subset\{1,...,N_E\},S_{i,e} \text{ occurs in Spring.}
\end{equation*}
\end{itemize}

\oldrev{Finally, these cost functions enable to formulate, after the single objective calibration problem \ref{inv_problem}, the following multi-objectives calibration problems.} The optimization problems taking into account signatures via the cost function defined in Eq. \ref{eq:totalcost} can be developed as Eq. \ref{eq:prl_singleoptim} for a single-objective optimization,
and as Eq. \ref{eq:prl_MOoptim-1} for a multi-objective optimization.
\begin{equation}\label{eq:prl_singleoptim}
\min_{\boldsymbol{\theta}\in\mathcal{O\subset\mathbb{R}}^{n}}\delta_{d}j_{d}(\boldsymbol{\theta}) +\delta_{c}\sum_{S_{i}\in N_{c}}\sigma_{S_{i}}\left|\frac{S_{i}^{s}(\boldsymbol{\theta})}{S_{i}^{o}(\boldsymbol{\theta})}-1\right| + \delta_{f}\sum_{S_{i}\in N_{f}}\sigma_{S_{i}}\frac{1}{N_E}\sum_{e=1}^{N_{E}}\left|\frac{S_{i,e}^{s}(\boldsymbol{\theta})}{S_{i,e}^{o}(\boldsymbol{\theta})}-1\right|
\end{equation}
\begin{equation}\label{eq:prl_MOoptim-1}
\min_{\boldsymbol{\theta}\in\mathcal{O\subset\mathbb{R}}^{n}} \Bigg( j_{d}(\boldsymbol{\theta}), \left\{\left|\frac{S_{i}^{s}(\boldsymbol{\theta})}{S_{i}^{o}(\boldsymbol{\theta})}-1\right|\right\}_{S_{i}\in N_{c}},
\left\{\frac{1}{N_E}\sum_{e=1}^{N_{E}}\left|\frac{S_{i,e}^{s}(\boldsymbol{\theta})}{S_{i,e}^{o}(\boldsymbol{\theta})}-1\right|\right\}_{S_{i}\in N_{f}}\Bigg)
\end{equation}
While the minimization problem with single-objective function \ref{eq:prl_singleoptim} is accessible for both global and distributed calibration methods, performing a multi-objective optimization as problem \ref{eq:prl_MOoptim-1} is sophisticated for distributed calibration considering a spatially distributed control vector adapted to a high dimensional hydrological optimization problems, and requiring a lot of cost gradient information.
\change{When using multi-objective optimization approaches in global calibration, it is possible to find a set of feasible solutions rather than a single optimal solution, as is the case in single-objective optimization (cf. Appendix \ref{MOO_appendix}). This is achieved through the use of a non-dominated sorting genetic algorithm (NSGA), which will be discussed in detail in Appendix \ref{GA_overview}. As a result, a set of non-inferior solutions, also known as a Pareto front, can be obtained (see Appendix \ref{pareto_front_appendix}). To select the optimal solution from the Pareto front, a method is proposed and depicted in Appendix \ref{select_sol_append}}.

Note that the objective functions $j_c$ and $j_f$ related to continuous and flood signatures have also been implemented in Fortran. This implementation imposes strict positivity of their components ($j^{S_i}_c$ and $j^{S_i}_f$) numerically to ensure that the total cost $J$ remains convex and differentiable. The numerical adjoint model has been also re-derived as needed by the variational calibration algorithm (refer to section \ref{subsec:forward_and_calibration_algos}). The cost function based flood event signatures $j_f$ can be computed thanks to a temporal mask of corresponding flood events selected by the segmentation algorithm, implemented in the Python \change{routines}, and passed to the Fortran \change{routines} via the wrapped interface.

\section{\oldrev{Data and numerical results analysis}\label{sec:Result}}

\oldrev{This section first presents the catchment-flood dataset used in this study. Next, flow signatures are analyzed via the comparison of observed and simulated signatures, in terms  \change{of} sensitivity \change{to model} parameters, and finally some are selected for signature-based model  calibration. The last part  analyzes the performances of model calibration with classical and signature-based metrics.
}

\subsection{\oldrev{Catchment information and data sources}}

\begin{figure}[ht!]
     \centering
     \includegraphics[scale=0.45]{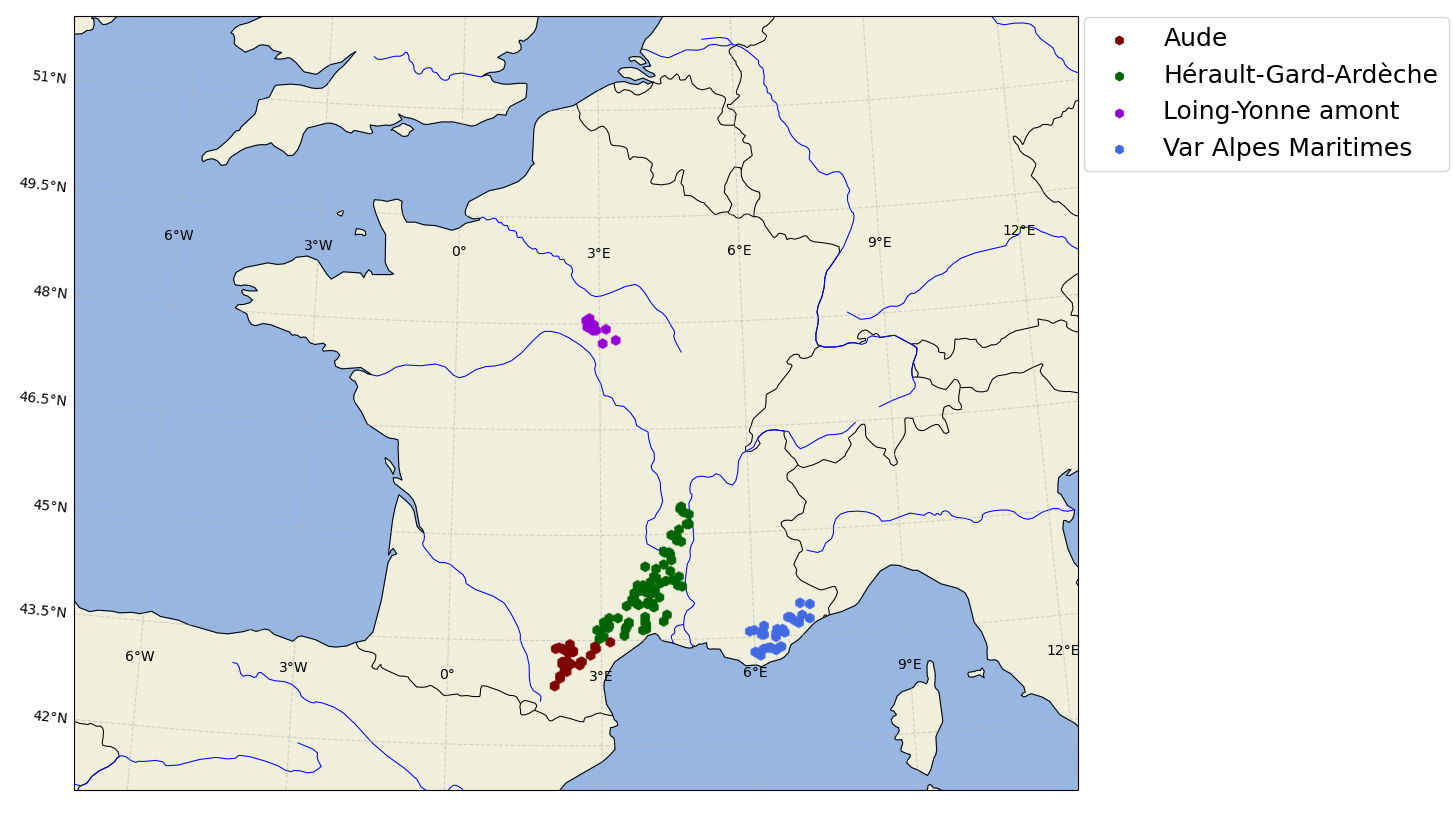}
     \caption{Spatial distribution of 141 catchments of the PICS dataset\oldrev{, consisting of 23 outlet gauges and 118 upstream gauges} on the map of France \oldrev{with four regions} denoted by different colors.}
     \label{fig:area_pics}
 \end{figure}

\oldrev{A relatively large dataset of catchment-floods mostly located in the French mediterranean region is used.}
This dataset \oldrev{stems from} \citet{jay2020estimation} \oldrev{and} contains time series of hydro-meteorological variables and time invariant catchment attributes for four high rainfall-flow areas in France, identified as study areas of the PICS research project\footnote{\href{https://pics.ifsttar.fr}{https://pics.ifsttar.fr}}. It encompasses 141 catchments including 23 outlet gauges, which are mostly located in the French Mediterranean region (Fig. \ref{fig:area_pics}). This is a subset of a larger dataset of 4,190 French catchments from INRAE-HYCAR research unit \citep{brigode2020data, delaigue2020data}. 
\oldrev{The hydrological model inputs consist of observation data, covering a period of about 13 years (2006 to 2019), that includes hourly distributed discharge and rainfall.} Discharge data  are  collected by the French Ministry of Environment covering the period of the forcing data and have been extracted from the (\href{http://www.hydro.eaufrance.fr/}{Hydro}) platform\footnote{\href{http://www.hydro.eaufrance.fr/}{http://www.hydro.eaufrance.fr/}}. The rainfall grids are the radar observation reanalysis ANTILOPE J+1 provided by Météo-France at a grid resolution of 1~km$^2$ \citep{champeaux2009mesures}. The potential evapotranspiration (PET) is obtained by applying a simple formula \citep{oudin2005potential} to SAFRAN\footnote{"Système d’Analyse Fournissant des Renseignements Atmosphériques à la Neige" in French} \citep{Quintana2008} temperature grids at 8~km resolution an empirically disaggregated at hourly time step and 1~km spatial resolution, i.e. at the same spatio-temporal resolution than rainfall. 
Note that observation data, rainfall grids and discharge time series, over the selected catchments have few missing data as detailed in Table \ref{tab:catch_info}, so that it can be neglected when performing the computations and analysis in this study. \oldrev{Table \ref{tab:catch_info} contains catchment information such as the river name, surface, code, number of upstream gauges, and missing rates in the outlet gauges.} Raster maps, at 1~km resolution, of upstream drained area and D8 flow directions have been obtained by processing fine DEM provided by IGN (Institut Geographique National).
\begin{table}[ht!] \renewcommand*{\arraystretch}{1.25}
\begin{centering}
\caption{General information about 23 outlet gauges of the PICS data. Code, river name, surface, missing rate of rainfall (\oldrev{respectively,} discharge) in outlet gauge during the period 2006-2019, and number of upstream are represented by the columns from left to right.}
\label{tab:catch_info}
\scalebox{0.8}{
\begin{tabular}{c|c|c|c|c}
    Code & River name & Surface ($km^2$) & Missing rates ($\%$) & Total upstream gauges\\
    \hline
    \hline
    H3201010 & Le Loing & 2302 & 0.14 (3.68) & 8 \\
    V3524010 & La Cance & 381 & 0.14 (4.31) & 3 \\
    V3744010 & Le Doux & 621 & 0.14 (4.02) & 2 \\
    V4154010 & L’Eyrieux & 649 & 0.14 (7.38) & 3 \\
    V5064010 & L’Ardèche & 2264 & 0.14 (4.22) & 9 \\
    V5474015 & La Cèze & 1112 & 0.14 (3.76) & 6 \\
    V7164015 & Le Gardon & 1093 & 0.14 (16.62) & 10 \\
    Y1232010 & L’Aude & 1828 & 0.14 (3.74) & 11 \\
    Y1364010 & Le Fresquel & 935 & 0.14 (3.74) & 4 \\
    Y1415020 & L’Orbiel & 242 & 0.14 (3.74) & 2 \\
    Y1564010 & L’Orbieu & 589 & 0.14 (3.77) & 3 \\
    Y1605050 & La Cesse & 251 & 0.14 (4.64) & 1 \\
    Y2332015 & L’Hérault & 2208 & 0.14 (7.22) & 12 \\ 
    Y2584010 & L’Orb & 1336 & 0.14 (4.04) & 11 \\
    Y3204040 & Le Lez & 168 & 0.14 (15.55) & 3 \\
    Y3444020 & Le Vidourle & 503 & 0.14 (7.97) & 4 \\
    Y3534010 & Le Vistre & 496 & 0.14 (4.42) & 1 \\
    Y4624010 & Le Gapeau & 535 & 0.14 (3.79) & 6 \\
    Y5312010 & L’Argens & 2512 & 0.14 (5.08) & 10 \\
    Y5444010 & La Giscle & 201 & 0.14 (9.96) & 2 \\
    Y5534030 & La Siagne & 492 & 0.14 (5.30) & 5 \\
    Y5615030 & Le Loup & 289 & 0.14 (3.79) & 1 \\
    Y6434010 & L’Estéron & 442 & 0.14 (7.70) & 1 \\
\end{tabular}}
\par\end{centering}
\end{table}

\subsection{\oldrev{Sensitivity analysis and selection of signatures for model calibration}}

To start with, \oldrev{the relative error is analyzed  between  observed signatures and simulated ones with a model calibrated using SBS algorithm and spatially uniform parameters.} Table \ref{tab:re_sign} shows that some hydrological signatures with a significant simulation error such as: $Cfp2$, $Cfp10$, $Cfp50$, $Elt$ and $Epf$ \oldrev{that could be better constrained with a signature-based calibration process as investigated in next subsection} (a list of all studied signatures \oldrev{with corresponding notations} is presented in Appendix \ref{appendix_sign}). 

Next, we survey the \oldrev{global} sensitivity of these signatures with the model parameters. We considered over 10,000 spatially uniform sets of the 6 model parameters, sampled using Saltelli generator \citep{saltelli2002making}, to estimate the \oldrev{total-order} Sobol indices across 23 gauged catchments (catchments downstream outlets of the dataset). \oldrev{Based on the results presented in Table \ref{tab:sens_sign}, it can be observed that the non conservative water exchange parameter $ml$ and the transfer parameter $c_{tr}$ exhibit the highest sensitivities to the studied signatures, both in terms of first-order and total-order. Our analysis suggests that these two parameters have the most significant impact on the output signatures as a result of their interactions with other inputs. This is in coherence with highest sensitivities found for soil depth and subsurface flow parameters of an event flash flood model found in \citet{garambois2013sobol, garambois2015characterization} on some catchments of the present set. Conversely, we found that parameters such as the interception $c_i$ and the production of runoff $c_p$ have little to no impact on the simulated signatures. We also observed that continuous signatures exhibit lower sensitivities than flood-event signatures in both first-order and total-order effects.} Furthermore, constraining hydrological model by flood event signatures along with a classical calibration metric (e.g. $1-NSE$ or $1-KGE$), \oldrev{which is based primarily on continuous records of streamflow}, is ideal to balance the model between the global score and the performance on flood events. We select for example \oldrev{the peak flow, denoted as $Epf$}, which is one of flood event signatures having both significant relative error and high sensitivity, to perform multi-critera calibration methods. Note that multi-criteria optimization methods with multi\oldrev{ple} signatures are absolutely reachable but will not be shown in this study for sake of brevity and simplify results analysis.

Furthermore, it is worth mentioning that global sensitivity analysis can be performed with local derivatives based approaches. A link between global Sobol indices and local derivatives has been proposed by \citet{sobol2010derivative} (refer also to \citet{lamboni2013derivative}). Global sensitivity matrices in three dimensions (sample size, parameters number, time) and sensitivity statistics, based on local derivatives computed by finite differences have been proposed in \citet{gupta2018revisiting, razavi2019multi} for geophysical models and applied to HBV-SASK lumped hydrologic model. Note that the VDA algorithm upgraded in the present work uses accurate local (in parameter space) cost function gradients, global in time and spatially distributed, computed with the adjoint method. Such method enables to compute accurate spatial sensitivity maps even for high dimensional parameter spaces (e.g. \citet{monnier2016inverse}) and deepening sensitivity analysis with our differentiable and spatially distributed hydrological model, along with accounting for sensitivity indices into the VDA algorithm, is a very interesting direction intentionally left for further research.

\begin{table}[ht!] \renewcommand*{\arraystretch}{1.25}
\begin{centering}
\caption{Relative error between simulated and observed signatures of the same model structure calibrated either with $1-NSE$ or $1-KGE$ by SBS algorithm for global optimization. The values (in the form of . [., .]) in each case represent respectively the median, mean and standard deviation of a signature over gauged catchments.}
\label{tab:re_sign}
\scalebox{0.8}{
\oldrev{
\begin{tabular}{cc|cc}
    \multirow{2}{*}{Notation} & \multirow{2}{*}{Signature type} & \multicolumn{2}{c}{Relative error on simulated signature} \\
    \cline{3-4}
    {} & {} & Cal. with $j_{d}^{NSE}$ & Cal. with $j_{d}^{KGE}$ \\
    \hline\hline
    Crc & \multirow{4}{*}{Continuous runoff coefficients} & 0.14 [0.28, 0.38] & 0.16 [0.3, 0.46] \\
    Crchf & {} & 0.24 [0.35, 0.35] & 0.26 [0.4, 0.45] \\
    Crclf & {} & 0.15 [0.3, 0.44] & 0.15 [0.33, 0.54] \\
    Crch2r & {} & 0.23 [0.4, 0.68] & 0.22 [0.38, 0.69] \\
    \hline
    Cfp2 & \multirow{4}{*}{Flow percentiles} & 0.72 [3.99, 21.14] & 0.76 [5.99, 29.98] \\
    Cfp10 & {} & 0.52 [2.64, 8.8] & 0.52 [2.87, 9.42] \\
    Cfp50 & {} & 0.29 [0.49, 0.85] & 0.2 [0.52, 0.99] \\
    Cfp90 & {} & 0.21 [0.37, 0.96] & 0.18 [0.38, 0.99] \\
    \hline
    Eff & Flood flow & 0.23 [0.32, 0.31] & 0.19 [0.31, 0.37] \\
    Ebf & Base flow & 0.22 [0.33, 0.39] & 0.22 [0.33, 0.41] \\
    \hline
    Erc & \multirow{4}{*}{Flood event runoff coefficients} & 0.2 [0.28, 0.26] & 0.18 [0.27, 0.26] \\
    Erchf & {} & 0.23 [0.32, 0.31] & 0.19 [0.31, 0.37] \\
    Erclf & {} & 0.22 [0.33, 0.39] & 0.22 [0.33, 0.41] \\
    Erch2r & {} & 0.12 [0.19, 0.2] & 0.13 [0.2, 0.24] \\
    \hline
    Elt & Lag time & 0.48 [0.96, 1.25] & 0.46 [0.82, 1.1] \\
    Epf & Peak flow & 0.28 [0.38, 0.35] & 0.25 [0.36, 0.41] \\
\end{tabular}}
}
\par\end{centering}
\end{table}

\begin{table}[ht!] \renewcommand*{\arraystretch}{1.25}
\begin{centering}
\caption{\oldrev{Median across gauged catchments of first- (respectively, total-) order variance-based sensitivity indices of the studied signatures to the model parameters.}}
\label{tab:sens_sign}
\scalebox{0.8}{
\oldrev{
\begin{tabular}{c|cccccc}
    \multirow{2}{*}{Signature} & \multicolumn{6}{c}{Model parameter} \\
    \cline{2-7}
    {} & $c_i$ & $c_p$ & $c_{tr}$ & $c_{tl}$ & $c_r$ & $ml$ \\
    \hline\hline
    Crc & -0.0 (0.0001) & -0.0004 (0.0004) & 0.1336 (1.2998) & 0.0006 (0.0002) & -0.0 (0.0) & 0.1167 (1.3778) \\
    Crchf & 0.0038 (0.0103) & 0.0268 (0.1155) & 0.3739 (0.8506) & 0.0153 (0.0123) & 0.0367 (0.1513) & 0.2245 (0.7919) \\
    Crclf & -0.0 (0.0001) & -0.0004 (0.0003) & 0.1299 (1.3018) & 0.0006 (0.0001) & -0.0 (0.0) & 0.1142 (1.3844) \\
    Crch2r & -0.0004 (0.0017) & 0.0193 (0.0255) & 0.1014 (0.1426) & 0.1099 (0.2055) & 0.1833 (0.2449) & 0.3984 (0.5481) \\
    Cfp2 & 0.0014 (0.056) & 0.0002 (0.0024) & 0.1283 (1.6008) & 0.0 (0.0) & 0.0 (0.001) & -0.0026 (1.2871) \\
    Cfp10 & -0.0 (0.0001) & -0.0002 (0.0001) & 0.128 (1.3353) & 0.0002 (0.0) & 0.0 (0.0) & 0.092 (1.3922) \\
    Cfp50 & -0.0001 (0.0001) & -0.0001 (0.0001) & 0.1267 (1.315) & 0.0005 (0.0001) & 0.0 (0.0) & 0.1043 (1.3933) \\
    Cfp90 & -0.0002 (0.0001) & -0.0006 (0.0015) & 0.1329 (1.2483) & 0.001 (0.0006) & -0.0001 (0.0002) & 0.1512 (1.3817) \\
    Eff & 0.0002 (0.0059) & 0.0699 (0.1939) & 0.306 (0.7807) & 0.0242 (0.022) & 0.0321 (0.1389) & 0.1872 (0.7303) \\
    Ebf & -0.0001 (0.001) & 0.0014 (0.0159) & 0.144 (1.1914) & 0.0018 (0.0019) & -0.0002 (0.0056) & 0.162 (1.3146) \\
    Erc & -0.0001 (0.0015) & 0.0076 (0.0314) & 0.18 (1.1633) & 0.0028 (0.0031) & -0.0001 (0.0011) & 0.1705 (1.2433) \\
    Erchf & 0.0002 (0.0059) & 0.0699 (0.1939) & 0.306 (0.7807) & 0.0242 (0.022) & 0.0321 (0.1389) & 0.1872 (0.7303) \\
    Erclf & -0.0001 (0.001) & 0.0014 (0.0159) & 0.144 (1.1914) & 0.0018 (0.0019) & -0.0002 (0.0056) & 0.162 (1.3146) \\
    Erch2r & 0.0057 (0.0099) & 0.0123 (0.0426) & 0.0873 (0.2124) & 0.0256 (0.0552) & 0.4387 (0.5797) & 0.1171 (0.2255) \\
    Elt & -0.0002 (0.0116) & -0.0004 (0.0293) & 0.0043 (0.087) & 0.0009 (0.0048) & 0.8832 (0.953) & 0.0127 (0.0568) \\
    Epf & -0.0008 (0.0026) & 0.0357 (0.1235) & 0.2505 (0.9199) & 0.0081 (0.0074) & 0.1099 (0.2632) & 0.1257 (0.8049) \\
\end{tabular}}}
\par\end{centering}
\end{table}

\subsection{\oldrev{Performance comparison of classical and signature-based calibration metrics}}

\change{In this section, we compare the performance of different models using both uniform and distributed optimization methods with different calibration metrics, including signature-based ones.} For spatially uniform calibration methods, we aim to compare different calibration metrics including classical single-objective optimization (CSOO), signature-based single-objective optimization (SSOO) and signature-based multi-objective optimization (SMOO). For \oldrev{spatially} distributed \oldrev{calibration} methods, two strategies selected for comparison are CSOO and SSOO.
In both \oldrev{spatially uniform or distributed calibration} scenarios, the models are calibrated on 23 outlet gauges of the PICS data on the calibration period 2006-2013. \oldrev{The validation of calibrated models performances is done in space and time following the three setups:}
\begin{itemize}
    \item on 23 outlet gauges on the validation period 2013-2019 (temporal validation, T\_Val),
    \item on 118 upstream gauges on the calibration period 2006-2013 (spatial validation, S\_Val),
    \item on 118 upstream gauges on the validation period 2013-2019 (spatio-temporal validation, S-T\_Val).
\end{itemize}

\subsubsection{\oldrev{Spatially uniform calibrations with NSGA}}
We first perform global calibrations using NSGA \change{with (i) classical single-objective functions; (ii) multi-criteria single-objective functions; and (iii) multi-objective functions}. Table \ref{tab:val_metric_uniform} displays the mean of different objective functions for calibration and validation (with 3 validation metrics), and for 3 optimization \change{methods} (CSOO, SSOO and SMOO) with various cost functions. In CSOO, we interpret that the model calibrated with $j_d^{KGE}=1-KGE$ produces a better result on the peak flow $j_f^{Epf}$, compared to the one calibrated with $j_d^{NSE}=1-NSE$. This explains why KGE criterion is more robust than NSE for constraining a hydrological model, since it is built on the decomposition of
NSE \citep{gupta2009decomposition}, which emphasizes relative importance of several hydrological features.
\oldrev{This finding is consistent with that of \citet{mizukami2019choice}. The authors calibrated daily models over numerous US catchments with multiple metrics, including NSE, weighted KGEs, annual peak flow bias (APFB), and they found that KGE resulted in better estimates of annual peak flows than NSE. Additionally, the best reproduction of annual peak flows was achieved with APFB, but this was at the expense of other high flow metrics.}

\begin{figure}[ht!]
\begin{centering}
\includegraphics[scale=0.85]{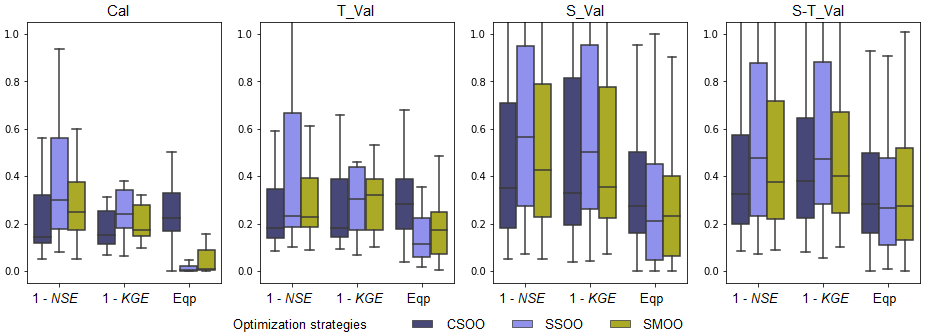}
\par\end{centering}
\caption{\oldrev{Comparison, with spatially uniform parameters, of calibration and validation metrics (optimal fit for $cost=0$) for three optimization approaches (CSOO, SSOO, SMOO) by constraining $1-NSE$ in case of global algorithms (SBS or NSGA). From left to right: calibration (Cal), temporal validation (T\_Val), spatial validation (S\_Val) and spatio-temporal validation (S-T\_Val).}}
\label{val_nse_uniform}
\end{figure}

\begin{figure}[ht!]
\begin{centering}
\includegraphics[scale=0.85]{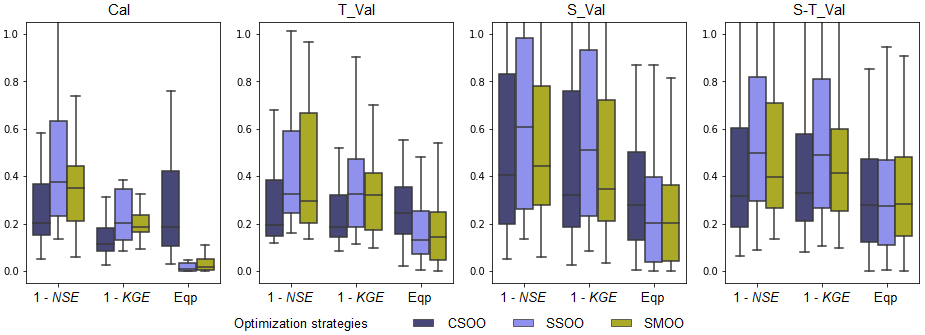}
\par\end{centering}
\caption{\oldrev{Comparison, with spatially uniform parameters, of calibration and validation metrics (optimal fit for $cost=0$) for three optimization approaches (CSOO, SSOO, SMOO) by constraining $1-KGE$ in case of global algorithms (SBS or NSGA)}. From left to right: calibration (Cal), temporal validation (T\_Val), spatial validation (S\_Val) and spatio-temporal validation (S-T\_Val).}
\label{val_kge_uniform}
\end{figure}

Using event signatures in addition to classical continuous  metrics in SSOO, we found that simulated peak flow is highly improved in terms of relative error (about $15$-$18$ times and $1.4$-$1.7$ times on average,  \oldrev{respectively,} for calibration and temporal validation) while classical calibrated metrics are significantly deteriorated (about $1.4$-$1.6$ times and $1.4$-$1.5$ times on average, \oldrev{respectively,} for calibration and temporal validation). This may arise from \oldrev{imbalances} between global score and performance in simulating flood event signature. \oldrev{To address this issue, careful consideration of the optimization weights assigned to objective functions is necessary in order to achieve a balance between model performance on short and long-term series. It should be noted that this approach can be time-consuming, as it requires numerous simulations to determine the appropriate optimization weights for the objective functions, typically using a L-curve approach. Alternatively, the use of global calibration algorithms, which do not require gradient information and can be solved using lower-dimensional optimization problems, can also address these imbalances through the application of a multi-objective optimization approach. This approach offers the advantage of keeping acceptable levels of deterioration of NSE and KGE while significantly improving the simulation of peak flow as shown by multi-objective SMOO results in} Fig. \ref{val_nse_uniform} and \ref{val_kge_uniform}. 

However, this global multi-\oldrev{objective} optimization algorithm is not capable to deal with high dimensional control vectors and the \oldrev{spatially} uniform \oldrev{parameter} setup here \oldrev{(under-parameterization)} led to unsatisfactory results in spatial and spatio-temporal validation metrics.
\oldrev{Therefore, a distributed calibration approach, such as using our VDA algorithm accounting for signatures, could improve the model performances. 
This approach maintains the same optimization weights as described above, and its performance will be evaluated in the subsequent section.}

\begin{table}[ht!] \renewcommand*{\arraystretch}{1.5}
\begin{centering}
\caption{Calibration, temporal, spatial and spatio-temporal validation metrics with spatially uniform calibrations with three strategies (CSOO, SSOO, SMOO) and global algorithms (SBS or NSGA), optimal fit for $cost=0$. The mean of calibration and validation cost values for different objective functions are displayed for each calibration metric \oldrev{- mean is computed over the catchment set: over the 23 outlet gauges for Cal and T\_Val, over the remaining 118 upstream gauges for S\_Val and S-T\_Val.}}
\label{tab:val_metric_uniform}
\scalebox{0.8}{
\begin{tabular}{c|c|cccc|cccc|cccc}
    \multirow{2}{*}{Method} & \multirow{2}{*}{Calibration metric} & \multicolumn{4}{c|}{$\overline{j_{d}^{NSE}}$} & \multicolumn{4}{c|}{$\overline{j_{d}^{KGE}}$} & \multicolumn{4}{c}{$\overline{j_{f}^{Epf}}$} \\
    \cline{3-14}
    {} & {} & Cal & T\_Val & S\_Val & S-T\_Val & Cal & T\_Val & S\_Val & S-T\_Val & Cal & T\_Val & S\_Val & S-T\_Val \\
    \hline\hline
    \rowcolor{lightgray}
    \cellcolor{white}
    \multirow{2}{*}{CSOO} & $j_{d}^{NSE}$ & 0.274 & 0.277 & 0.901 & 0.616 & 0.239 & 0.369 & 0.687 & 0.736 & 0.279 & 0.324 & 0.387 & 0.357 \\
     & $j_{d}^{KGE}$ & 0.352 & 0.330 & 1.048 & 0.795 & 0.183 & 0.323 & 0.665 & 0.721 & 0.267 & 0.280 & 0.379 & 0.344 \\
    \hline
    \rowcolor{lightgray}
    \cellcolor{white}
    \multirow{2}{*}{SSOO} & $j_{d}^{NSE}/2+j_{f}^{Epf}/2$ & 0.447 & 0.418 & 1.056 & 0.889 & 0.377 & 0.476 & 0.759 & 0.853 & 0.014 & 0.189 & 0.346 & 0.372 \\
     & $j_{d}^{KGE}/2+j_{f}^{Epf}/2$ & 0.551 & 0.431 & 1.259 & 0.956 & 0.335 & 0.443 & 0.777 & 0.833 & 0.017 & 0.209 & 0.337 & 0.358 \\
    \hline
    \rowcolor{lightgray}
    \cellcolor{white}
    \multirow{2}{*}{SMOO} & $\{j_{d}^{NSE},j_{f}^{Epf}\}$ & 0.341 & 0.351 & 1.020 & 0.845 & 0.271 & 0.420 & 0.703 & 0.803 & 0.087 & 0.215 & 0.336 & 0.391 \\
     & $\{j_{d}^{KGE},j_{f}^{Epf}\}$ & 0.456 & 0.409 & 1.163 & 0.821 & 0.243 & 0.368 & 0.683 & 0.724 & 0.048 & 0.182 & 0.316 & 0.389 \\
\end{tabular}}
\par\end{centering}
\end{table}

\begin{figure}[ht!]
\begin{centering}
\includegraphics[scale=0.77]{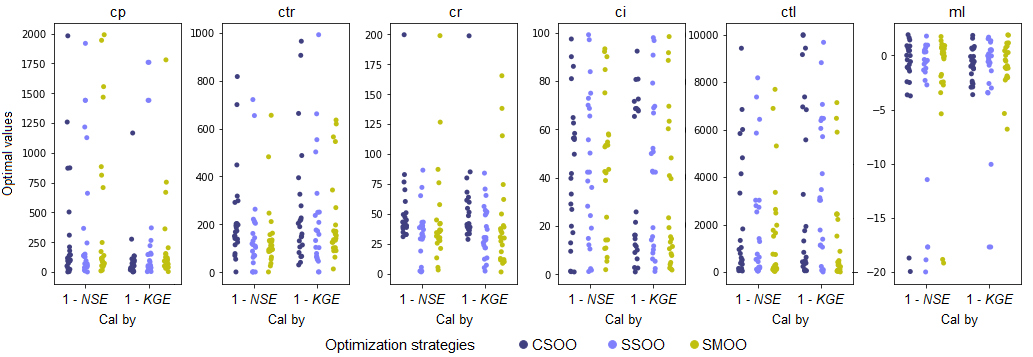}
\par\end{centering}
\caption{\oldrev{Analysis of spatially uniform calibrated parameters over the whole catchment sample. In each scatterplot,} the first column present the distribution of a parameter for 3 optimization strategies (CSOO, SSOO, SMOO) using $j_d^{NSE}$, whereas the strategies in the second column use $j_d^{KGE}$ as the dominant (or constrained) objective function. The boundary conditions of the model parameters are given in Appendix \ref{BC_6param}.}
\label{param_space_uniform}
\end{figure}

\oldrev{As shown in Fig. \ref{param_space_uniform}, the corresponding optimal parameters obtained using various optimization strategies are presented. Based on our preliminary analysis, it is evident that the distribution over studied catchments of $c_r$ has an important difference when performing traditional calibration (CSOO) and multi-criteria calibration methods (SSOO and SMOO).} We recall that $c_r$ is the routing parameter in our conceptual design (Fig. \ref{fig:SMASH_grid_model}), so it has a crucial role in producing the peak flow $Epf$. Additionally, the sensitivity analysis in Table \ref{tab:sens_sign} has indicated that $c_r$ is one of the three parameters explaining most of the sensitivity of the peak flow. 

\oldrev{The above result on the importance of lateral flow components in a flood hydrological model is in coherence with existing works, for example as shown in \citet{garambois2013sobol} on few catchments-flood events used in the present study, in addition to high sensitivity to subsurface flow parameter (see also \citet{douinot2018using}) the temporal sensitivity of kinematic wave compound friction parameters in a distributed flash flood model increases with flood magnitude. Improving hydraulic  meaningfulness of hydrological models is an important topic since it can improve floods discharge modeling in high resolution catchment-flood models (e.g. \citet{bout2018validity, li2021crest, kirstetter2021b} with shallow water models and simplifications) but also improve internal state-flux coherence and realism as required for instance to assimilate remote sensing observables of river suface such as height and width (e.g. \citet{paiva2011large, pujol2020estimation, pujol2022multi}).}

\subsubsection{\oldrev{Spatially distributed calibrations with VDA algorithm}}

\begin{figure}[ht!]
\begin{centering}
\includegraphics[scale=0.85]{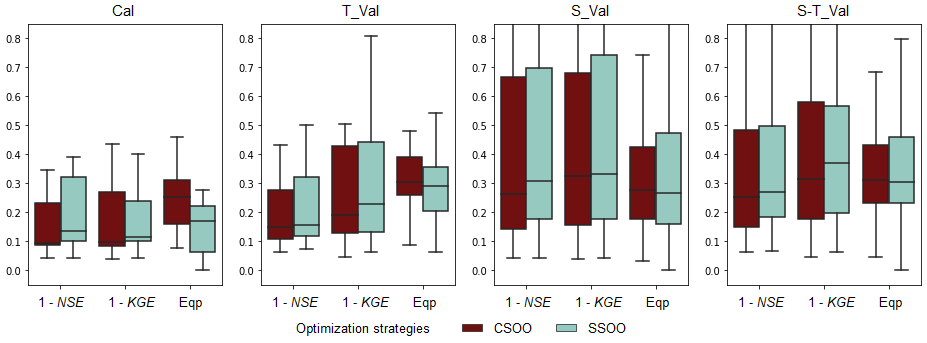}
\par\end{centering}
\caption{\oldrev{Comparison, with spatially distributed parameters, of calibration and validation metrics (optimal fit for $cost=0$) for two optimization approaches (CSOO, SSOO) by constraining $1-NSE$ in case of distributed calibration. From left to right: calibration (Cal), temporal validation (T\_Val), spatial validation (S\_Val) and spatio-temporal validation (S-T\_Val).}}
\label{val_nse_dis}
\end{figure}

\begin{figure}[ht!]
\begin{centering}
\includegraphics[scale=0.85]{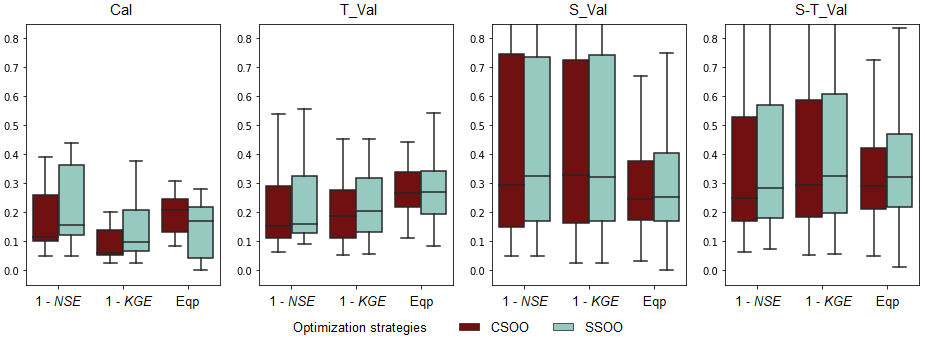}
\par\end{centering}
\caption{\oldrev{Comparison, with spatially distributed parameters, of calibration and validation metrics (optimal fit for $cost=0$) for two optimization approaches (CSOO, SSOO) by constraining $1-KGE$ in case of distributed calibration. From left to right: calibration (Cal), temporal validation (T\_Val), spatial validation (S\_Val) and spatio-temporal validation (S-T\_Val).}}
\label{val_kge_dis}
\end{figure}

Now, spatially distributed calibrations with the VDA algorithm using multi-criteria cost \change{functions, including signatures,} are performed. We employ SSOO technique for a distributed calibration using L-BFGS-B algorithm provided a first guess by SBS algorithm. In overall, all of obtained scores in Table \ref{tab:val_metric_distributed} are significantly enhanced compared to the uniform calibration method, thanks to spatially distributed control vectors granting more flexibility to reproduce observed discharge. Instead of a sharp decline of $j_f^{Epf}$ as above, this relative error slightly decreases about $1.5$ times (from about $0.25$ down to $0.16$) in calibration and from about $0.32$ down to $0.28$ in temporal validation, but instead, the scores (NSE and KGE) are slightly reduced in calibration and have an inappreciable deterioration in temporal validation. So in this case, we do not have \oldrev{imbalances} between the model performances on short and long-term series when employing SSOO. We observe clearly in Fig. \ref{val_nse_dis} and \ref{val_kge_dis} that the error of simulated pick flow is significantly reduced while the deterioration level of the scores remains tolerable, particularly in calibration and temporal validation.
\begin{table}[ht!] \renewcommand*{\arraystretch}{1.5}
\begin{centering}
\caption{Calibration, temporal, spatial and spatio-temporal validation metrics with spatially distributed control vectors (optimal fit for $cost=0$). The mean of calibration and validation cost values for different objective functions are displayed for each calibration metric.}
\label{tab:val_metric_distributed}
\scalebox{0.8}{
\begin{tabular}{c|c|cccc|cccc|cccc}
    \multirow{2}{*}{Method} & \multirow{2}{*}{Calibration metric} & \multicolumn{4}{c|}{$\overline{j_{d}^{NSE}}$} & \multicolumn{4}{c|}{$\overline{j_{d}^{KGE}}$} & \multicolumn{4}{c}{$\overline{j_{f}^{Epf}}$} \\
    \cline{3-14}
    {} & {} & Cal & T\_Val & S\_Val & S-T\_Val & Cal & T\_Val & S\_Val & S-T\_Val & Cal & T\_Val & S\_Val & S-T\_Val \\
    \hline\hline
    \rowcolor{lightgray}
    \cellcolor{white}
    \multirow{2}{*}{CSOO} & $j_{d}^{NSE}$ & 0.221 & 0.244 & 0.655 & 0.596 & 0.233 & 0.355 & 0.553 & 0.597 & 0.274 & 0.334 & 0.381 & 0.376 \\
     & $j_{d}^{KGE}$ & 0.239 & 0.231 & 0.802 & 0.702 & 0.140 & 0.292 & 0.617 & 0.701 & 0.226 & 0.295 & 0.365 & 0.364 \\
    \hline
    \rowcolor{lightgray}
    \cellcolor{white}
    \multirow{2}{*}{SSOO} & $j_{d}^{NSE}/2+j_{f}^{Epf}/2$ & 0.251 & 0.241 & 0.831 & 0.639 & 0.231 & 0.305 & 0.586  & 0.612 & 0.183 & 0.298 & 0.392 & 0.383 \\
     & $j_{d}^{KGE}/2+j_{f}^{Epf}/2$ & 0.297 & 0.245 & 0.964 & 0.671 & 0.190 & 0.300 & 0.617 & 0.647 & 0.152 & 0.271 & 0.376 & 0.387 \\
\end{tabular}}
\par\end{centering}
\end{table}

\begin{figure}[ht!]
\begin{centering}
\includegraphics[scale=0.85]{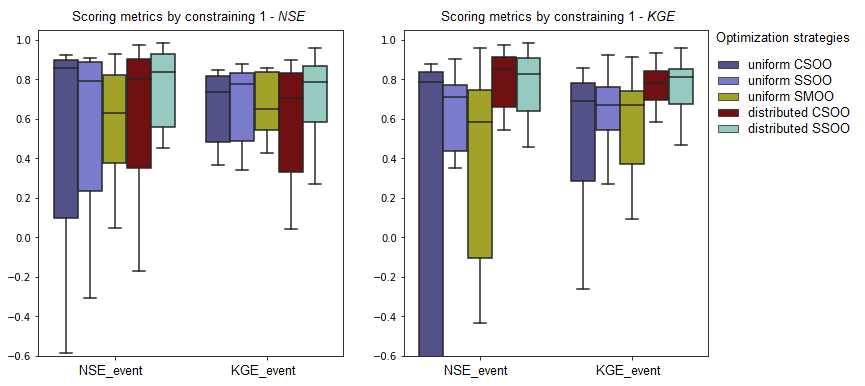}
\par\end{centering}
\caption{Comparison of scoring metrics computed on 111 events picked from 23 outlet gauges on the \oldrev{calibration} period 2006-2013 for \oldrev{the five optimization strategies}, by constraining $1-NSE$ (left) and $1-KGE$ (right), optimal fit for $score=1$.}
\label{cost_event}
\end{figure}

Ultimately, the scoring metrics are computed on 111 flood events picked from 23 outlet gauges (by segmentation method depicted in Algorithm \ref{algo_seg}) on the calibration period. The results plotted in Fig. \ref{cost_event} show that, in distributed calibration, the score of constrained calibration metric is not decreased but even improved from 0.80 (\oldrev{respectively,} 0.71) up to 0.83 (\oldrev{respectively,} 0.78) in median for NSE (\oldrev{respectively,} KGE). It indicates that the optimum of the model parameters has moved to another location that produces a better performance in simulating flood events by slightly reducing the scores in simulating the low-flow.

\begin{figure}[ht!]
\begin{centering}
\includegraphics[scale=0.155]{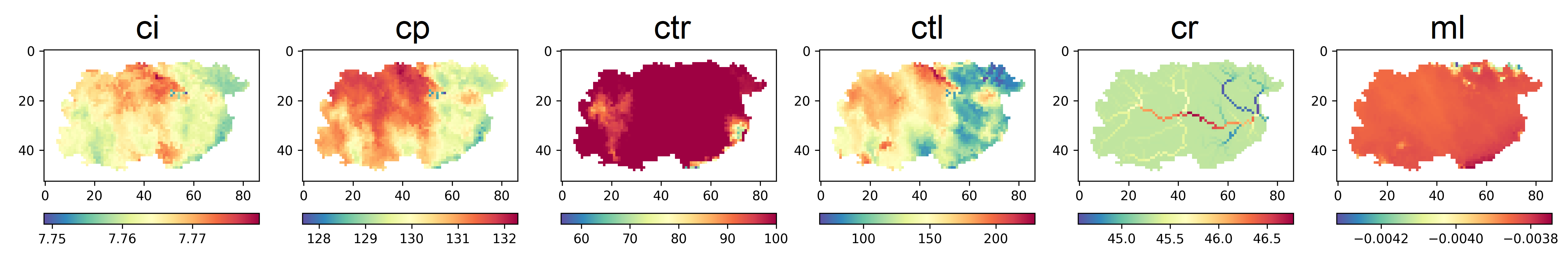}
\par\end{centering}
\caption{\oldrev{Spatially distributed optimal parameters ($\boldsymbol{\hat{\theta}}(x)\equiv\left(c_{i}(x),c_{p}(x),c_{tr}(x),c_{tl}(x),c_{r}(x),ml(x)\right)^{T}$) for the Argens River basin, obtained by minimizing $j_{d}^{KGE}/2 + j_{f}^{Epf}/2$.}}
\label{fig:map_Y5312010}
\end{figure}

Regarding to the parameter space, Table \ref{tab:first_guess} presents \change{statistical quantities, including mean, median and standard deviation,} of the spatially uniform parameter sets obtained using 4 calibration metrics for the studied catchments. Comparing to the spatially distributed optimal parameters in Table \ref{tab:dis_space_param}, we interpret that the mean of distributed parameters over all catchments in the 4 cases (corresponding to 4 calibration metrics) is globally coherent to the distribution of the first guess. Several parameters are almost spatially uniform (e.g. the non conservative water exchange parameter $ml$ has a small distributed deviation in median (\oldrev{respectively,} in average) over all catchments $0.01$ (\oldrev{respectively,} $0.05$) (calibrated with $j_d^{NSE}$) compared to its distributed mean in median (\oldrev{respectively,} in average) $-0.59$ (\oldrev{respectively,} $-4.98$)). Conversely, the transfer parameter $c_{tl}$ has a great distributed deviation (in median over all catchments) $193.72$ compared to its distributed average $347.87$, that also has a massive difference to its distributed median $114.79$.
Fig. \ref{fig:map_Y5312010} illustrates the spatially distributed optimal parameters at \oldrev{the largest catchment (the Argens River)}, for a distributed calibration with $j_{d}^{KGE}/2 + j_{f}^{Epf}/2$.

\oldrev{Reducing the over-parameterization in distributed hydrological models calibration problems through spatial constrains while enhancing regional parameters consistency remains a key issue, especially for flash flood prediction at ungauged locations (e.g. classical post-regionalization in \citet{garambois2015characterization} on French Mediterranean flash floods). 
This issue can be tackled with calibration approaches accounting for physiographic descriptors through 
regularizations (e.g. \citet{de2019regularization, jay2022spatially} in multi-gauges calibration problems) or through pre-regionalization mappings, such as the multi-scale parameter regionalization approach (MPR) from \citet{samaniego2010multiscale}, used for example in \citet{mizukami2017towards}. In addition to exploiting the information of multi-scale signatures in calibration with the present VDA algorithm, the use of a pre-regionalization scheme, i.e. "strong constrains" in the forward model in form of a mapping between physiographic covariables and conceptual hydrological parameter fields, represent an interesting perspective for future research.}

\begin{table}[ht!] \renewcommand*{\arraystretch}{1.5}
\begin{centering}
\caption{Uniform optimal parameters calibrated by SBS algorithm with 4 calibration metrics for each catchment on its outlet gauge. The values (in the form of . [., .]) in each case
represent respectively the median, mean and standard deviation of the optimal parameter values over all catchments of the dataset.}
\label{tab:first_guess}
\scalebox{0.8}{
\begin{tabular}{c|cc|cc}
    \multirow{2}{*}{Parameter} & \multicolumn{4}{c}{Calibration metric} \\
    \cline{2-5}
    {} & $j_{d}^{NSE}$ & $j_{d}^{NSE}/2 + j_{f}^{Epf}/2$ & $j_{d}^{KGE}$ & $j_{d}^{KGE}/2 + j_{f}^{Epf}/2$ \\
    \hline\hline
    $c_i$  & 14.71 [20.3, 26.22]  & 16.93 [20.83, 26.48]  & 17.6 [27.15, 33.07]  & 17.27 [30.17, 35.83]  \\
    $c_p$  & 169.99 [291.17, 434.58]  & 146.04 [310.14, 505.68]  & 151.87 [286.79, 483.33]  & 141.56 [289.69, 466.99]  \\
    $c_{tr}$  & 171.76 [286.6, 269.5]  & 162.66 [313.49, 304.83]  & 266.32 [431.2, 355.04]  & 267.21 [436.83, 360.62]  \\
    $c_{tl}$  & 347.87 [812.15, 1274.12]  & 250.42 [1366.73, 2789.22]  & 383.51 [1413.93, 2749.35]  & 262.89 [1337.96, 2795.16]  \\
    $c_r$  & 41.32 [52.63, 34.05]  & 40.94 [50.97, 30.97]  & 41.33 [51.2, 30.29]  & 40.24 [50.2, 27.58]  \\
    $ml$  & -0.59 [-4.98, 8.21]  & 0.0 [-3.81, 7.34]  & -0.0 [-3.62, 7.31]  & -0.0 [-3.28, 6.39]  \\
\end{tabular}}
\par\end{centering}
\end{table}
\begin{table}[ht!] \renewcommand*{\arraystretch}{1.5}
\begin{centering}
\caption{Analysis of spatially distributed parameter sets of the models corresponding to 4 calibration metrics. First, spatial median ($\tilde{.}$), average ($\overline{.}$) and standard deviation ($\sigma_{.}$) for each parameter field are calculated for each catchment, then their median, mean and standard deviation over all catchments are represented in the form of . [., .].}
\label{tab:dis_space_param}
\scalebox{0.8}{
\begin{tabular}{c|cc|cc}
    \multirow{2}{*}{Parameter} & \multicolumn{4}{c}{Calibration metric} \\
    \cline{2-5}
    {} & $j_{d}^{NSE}$ & $j_{d}^{NSE}/2+j_{f}^{Epf}/2$ & $j_{d}^{KGE}$ & $j_{d}^{KGE}/2+j_{f}^{Epf}/2$ \\
    \hline\hline
    $\tilde{c_i}$  & 15.45 [20.22, 26.34] & 10.91 [20.14, 26.75] & 17.6 [27.19, 33.16] & 17.3 [30.32, 36.05] \\
    $\overline{c_i}$  & 14.71 [20.3, 26.22] & 16.93 [20.83, 26.48] & 17.6 [27.15, 33.07] & 17.27 [30.17, 35.83] \\
    $\sigma_{c_i}$  & 0.22 [0.91, 1.46] & 0.07 [1.13, 4.16] & 0.13 [0.57, 1.28] & 0.05 [0.82, 3.09] \\
    \hline
    $\tilde{c_p}$  & 161.81 [286.05, 435.48] & 145.79 [314.19, 518.0] & 156.65 [288.75, 476.53] & 148.27 [296.94, 485.79] \\
    $\overline{c_p}$ & 169.99 [291.17, 434.58] & 146.04 [310.14, 505.68] & 151.87 [286.79, 483.33] & 141.56 [289.69, 466.99] \\
    $\sigma_{c_p}$ & 38.52 [60.58, 59.64] & 8.95 [37.57, 44.15] & 31.08 [53.49, 57.35] & 12.03 [46.82, 102.47] \\
    \hline
    $\tilde{c_{tr}}$  & 174.6 [287.44, 270.08] & 158.48 [317.5, 311.03] & 266.09 [429.0, 353.73] & 267.12 [447.27, 372.32] \\
    $\overline{c_{tr}}$  & 171.76 [286.6, 269.5] & 162.66 [313.49, 304.83] & 266.32 [431.2, 355.04] & 267.21 [436.83, 360.62] \\
    $\sigma_{c_{tr}}$  & 13.88 [28.84, 35.82] & 3.23 [25.3, 57.2] & 5.68 [24.28, 60.74] & 1.45 [24.39, 74.56] \\
    \hline
    $\tilde{c_{tl}}$  & 114.79 [675.7, 1276.32] & 127.09 [1322.22, 2806.59] & 180.63 [1332.45, 2784.42] & 146.96 [1322.17, 2803.08] \\
    $\overline{c_{tl}}$  & 347.87 [812.15, 1274.12] & 250.42 [1366.73, 2789.22] & 383.51 [1413.93, 2749.35] & 262.89 [1337.96, 2795.16] \\
    $\sigma_{c_{tl}}$  & 193.72 [355.92, 433.62] & 34.67 [139.21, 252.41] & 69.91 [222.54, 388.17] & 31.82 [61.87, 81.06] \\
    \hline
    $\tilde{c_r}$  & 41.37 [52.08, 34.42] & 41.37 [52.08, 34.42] & 41.37 [52.08, 34.42] & 41.37 [52.08, 34.42] \\
    $\overline{c_r}$  & 41.32 [52.63, 34.05] & 40.94 [50.97, 30.97] & 41.33 [51.2, 30.29] & 40.24 [50.2, 27.58] \\
    $\sigma_{c_r}$  & 4.66 [6.01, 5.04] & 1.45 [5.17, 10.17] & 3.01 [5.31, 10.29] & 1.34 [5.08, 13.82] \\
    \hline
    $\tilde{ml}$  & -0.59 [-4.98, 8.21] & 0.0 [-3.72, 7.42] & 0.0 [-3.61, 7.31] & -0.0 [-3.27, 6.39] \\
    $\overline{ml}$  & -0.59 [-4.98, 8.21] & 0.0 [-3.81, 7.34] & -0.0 [-3.62, 7.31] & -0.0 [-3.28, 6.39] \\
    $\sigma_{ml}$  & 0.01 [0.05, 0.09] & 0.0 [0.17, 0.73] & 0.02 [0.05, 0.09] & 0.0 [0.08, 0.34] \\
\end{tabular}}
\par\end{centering}
\end{table}

\section{Conclusion\label{sec:Conclusion}}

\oldrev{In this study, we enhanced the calibration process of the conceptual distributed hydrological model SMASH for Mediterranean floods by incorporating hydrological signatures and various multi-criteria optimization strategies. First, we computed and analyzed both continuous signatures and flood event signatures. Subsequently, we used sensitivity analysis to select appropriate signatures for constraining the model. Finally, we performed signatures-based multi-criteria optimization approaches, which demonstrated their robustness and reliability in improving simulated peak flood events without significantly compromising the NSE and KGE. Notably, for distributed calibration, the model constrained by the signature performed better in simulating flood events and achieved higher NSE and KGE scores compared to the model calibrated without using signatures. These results highlight the superiority of signature-based calibration approaches, particularly in flash flood prediction. Furthermore, we compared the parameter spaces of different models to provide insights into the optimal transition from traditional calibration approaches to signature-based calibration methods.}

Our proposed calibration strategy addresses the need for an intelligent approach to model calibration in the presence of multiple objectives and complex hydrological processes. This approach offers a new perspective on \change{the hydrological calibration process, accounting for the incorporation of classical discharge metrics and multi-scale signatures, that can provide a more comprehensive assessment of the model performance. As a potential perspective of improvement, our method} could be reinforced via the use of multi-source information such as from remotely sensed data products and of multi-gauge streamflow series in regionalization problems. The segmentation algorithm could be tested on larger flood samples, also including catchment rainfall moments \citep{zoccatelli2011spatial, emmanuel2015influence} describing rainfall patterns for floods analysis (e.g. \citet{garambois2014analysis, saharia2021impact}) and in order to prepare learning sets for training hybrid flood modeling-correction approaches. \change{Building on these insights, future work will aim to address the issue of equifinality, previously mentioned in the introduction section. The VDA algorithm could be upgraded with Bayesian elements, accounting for signatures and sensitivity. We will also develop better spatial constraints by mapping physiographic descriptors to parameters fields (pre-regionalization). The effectiveness of new approaches could be evaluated and analyzed over large datasets with signatures thanks to our segmentation algorithm. The proposed method could be easily applied to new spatially distributed model structures and hypothesis testing.}

\appendix

\section{Classical calibration metrics in hydrology}\label{append:metrics}

Nash--Sutcliffe efficiency (NSE):
\begin{equation}\label{eq:appendix_NSE}
NSE=1-\frac{\sum_{t=0}^{T}(Q(t)-Q^{*}(t))^{2}}{\sum_{t=0}^{T}(Q^{*}(t)-\overline{Q^{*}})^{2}}
\end{equation}
where $Q(t)$ is the simulated discharge at time
$t$, $Q^{*}(t)$ is the observed discharge at time $t$ and $\overline{Q^{*}}$
is the mean observed discharge.

Kling--Gupta efficiency (KGE):
\begin{equation}\label{eq:appendix_KGE}
KGE=1-\sqrt{\alpha(r-1)^{2}+\beta(\frac{\sigma}{\sigma^{*}}-1)^{2}+\gamma(\frac{\mu}{\mu^{*}}-1)^{2}}
\end{equation}
where $r$ is the linear correlation between observations
and simulations, $\sigma$ and $\sigma^{*}$ are the standard
deviation in simulations and observations, respectively, $\mu$
and $\mu^{*}$ are the mean discharge in simulations and observations,
respectively, and $\alpha,\beta,\gamma$ are the optimization weight parameters.
\section{List of studied signatures}    
\label{appendix_sign}
Denote $P(t)$ and $Q(t)$ are the rainfall and runoff at time $t\in\mathbf{U}$, where $\mathbf{U}$ is the study period. Then $Qb(t)$ and $Qq(t)$ are the baseflow and quickflow computed using a classical technique for streamflow separation (please refer to \citet{lyne1979stochastic} and \citet{nathan1990evaluation} for more details). 
The continuous signatures are calculated \oldrev{over the entire study period} as Table \ref{tab:continuous_sign}. 
For an event that occurs within a period $\mathbf{E} \subset \mathbf{U}$, the flood event signatures are calculated as Table \ref{tab:event_sign}.

\begin{table}[ht!] \renewcommand*{\arraystretch}{1.25}
\begin{centering}
\caption{List of all studied continuous signatures.}
\label{tab:continuous_sign}
\scalebox{0.8}{
\begin{tabular}{c c m{8cm} c c}
    Notation & Signature & Description & Formula & Unit\\
    \hline\hline
    Crc & \multirow{4}{*}{Runoff coefficients} & Coefficient relating the amount of runoff to the amount of precipitation received & $\frac{\int^{t\in\mathbf{U}} Q(t)dt}{\int^{t\in\mathbf{U}} P(t)dt}$ & -\\
    Crchf & {} & Coefficient relating the amount of high-flow to the amount of precipitation received & $\frac{\int^{t\in\mathbf{U}} Qq(t)dt}{\int^{t\in\mathbf{U}} P(t)dt}$ & -\\
    Crclf & {} & Coefficient relating the amount of low-flow to the amount of precipitation received & $\frac{\int^{t\in\mathbf{U}} Qb(t)dt}{\int^{t\in\mathbf{U}} P(t)dt}$ & -\\
    Crch2r & {} & Coefficient relating the amount of high-flow to the amount of runoff & $\frac{\int^{t\in\mathbf{U}} Qq(t)dt}{\int^{t\in\mathbf{U}} Q(t)dt}$ & -\\
    \hline
    Cfp2 & \multirow{4}{*}{Flow percentiles} & \multirow{4}{*}{$2\%$, $10\%$, $50\%$ and $90\%$-quantiles from flow duration curve} & $quantile(Q(t), 0.02)$ & \multirow{4}{*}{$mm$}\\
    Cfp10 & {} & {} & $quantile(Q(t), 0.1)$ & {}\\
    Cfp50 & {} & {} & $quantile(Q(t), 0.5)$ & {}\\
    Cfp90 & {} & {} & $quantile(Q(t), 0.9)$ & {}\\
\end{tabular}}
\par\end{centering}
\end{table}

\begin{table}[ht!] \renewcommand*{\arraystretch}{1.25}
\begin{centering}
\caption{List of all studied flood event signatures.}
\label{tab:event_sign}
\scalebox{0.8}{
\begin{tabular}{c c m{8cm} c c}
    Notation & Signature & Description & Formula & Unit\\
    \hline\hline
    Eff & Flood flow & Amount of quickflow in flood event & $\int^{t\in\mathbf{E}} Qq(t)dt$ & $mm$\\
    \hline
    Ebf & Base flow & Amount of baseflow in flood event & $\int^{t\in\mathbf{E}} Qb(t)dt$ & $mm$\\
    \hline
    Erc & \multirow{4}{*}{Runoff coefficients} & Coefficient relating the amount of runoff to the amount of precipitation received & $\frac{\int^{t\in\mathbf{E}} Q(t)dt}{\int^{t\in\mathbf{E}} P(t)dt}$ & -\\
    Erchf & {} & Coefficient relating the amount of high-flow to the amount of precipitation received & $\frac{\int^{t\in\mathbf{E}} Qq(t)dt}{\int^{t\in\mathbf{E}} P(t)dt}$ & -\\
    Erclf & {} & Coefficient relating the amount of low-flow to the amount of precipitation received & $\frac{\int^{t\in\mathbf{E}} Qb(t)dt}{\int^{t\in\mathbf{E}} P(t)dt}$ & -\\
    Erch2r & {} & Coefficient relating the amount of high-flow to the amount of runoff & $\frac{\int^{t\in\mathbf{E}} Qq(t)dt}{\int^{t\in\mathbf{E}} Q(t)dt}$ & -\\
    \hline
    Elt & Lag time & Difference time between the peak runoff and the peak rainfall & $\arg\max_{t\in\mathbf{E}} Q(t) - \arg\max_{t\in\mathbf{E}} P(t)$ & $h$\\
    \hline
    Epf & Peak flow & Peak runoff in flood event & $\max_{t\in\mathbf{E}} Q(t)$ & $mm$\\
\end{tabular}}
\par\end{centering}
\end{table}

\section{Segmentation algorithm}\label{app:sec:seg-algo}
\change{The proposed segmentation algorithm is illustrated in Algorithm \ref{algo_seg}. First, we identify event peak discharges using a peak detection algorithm, which allows for several parameters to be set, such as minimum peak height (mph) or minimum distance between two successive peaks (mpd), among others \citep{marcos_duarte_2021_4599319}. For instance, we consider events that exceed the 0.995-quantile of the discharge as important events (mph criterion), and events are considered to be distinct if they are separated by at least 12 hours (mpd criterion). 
Subsequently, we determine the starting and ending dates for each event. The starting date of the event is considered to be the moment when the rain starts to increase dramatically, which is sometime 72 hours before the peak discharge. To calculate this, we compute the gradient of the rainfall and choose the peaks of rainfall gradient that exceed the 0.8-quantile. These peaks correspond to the moments when there is a sharp increase in rainfall. However, we also require an additional criterion called the "energy criterion", which takes into account the "rainfall energy" for a more robust detection of flood start time. The rainfall energy is computed as the sum of squares of the rainfall observed in a 24-hour period, counted from 1 hour before the peak of rainfall gradient. The starting date is the first moment when the rainfall energy exceeds 0.2 of the maximal rainfall energy observed in the 72-hour period before the peak discharge, based on the gradient criterion. Finally, we aim to find the ending date by using baseflow separation. We
compute the difference between the discharge and its baseflow from
the peak discharge until the end of \oldrev{study period} (which lasts for 10 days from the starting date of the event). The ending date is the moment when the difference between the discharge and its baseflow is minimal in a 48-hour period, counted from 1 hour before this moment. Note that these values are adapted to the basins and flood scales studied.}
\begin{algorithm}
\change{\begin{flushleft}
For each catchment, considering 2 time series $(T,Q)$ and $(T,P)$ where:\\ 
$T=(t_{1},...,t_{n})$ is time (by hour),
$Q=(q_{1},...,q_{n})$ is the discharge, and
$P=(p_{1},...,p_{n})$ is the rainfall.
\end{flushleft}
\begin{enumerate}
\item Detecting peaks that exceed the 0.995-quantile of the discharge, that
can be considered as important events:\\ 
$E=(t_{i})_{1\leq i\leq n}$ s.t. $q_{i}>Quant_{0.995}(Q)$
\item For each event $t_{j}\in E$:
\begin{enumerate}
\item Determining a starting date based on the ``rainfall gradient criterion''
and the ``rainfall energy criterion'':
\begin{enumerate}
\item Selecting rainfalls gradient those exceed its $0.8$-quantile, considered
as the ``rainfall events'':\\
$RE=(t_{k})_{t_{k}\in(t_{j}-72,t_{j})}$ 
s.t. $\nabla P(t_{k})>Quant_{0.8}(\nabla P([t_{j}-72,t_{j}]))$
\item Defining the rainfall energy function:\\
$f(t_{x})=||(p_{x}-1,...,p_{x}+23)||_{2}$\\
then the starting date is the first moment the rainfall energy exceeds
0.2 of the maximal rainfall energy:\\
$sd=\min(t_{s})_{t_{s}\in RE}$ 
s.t. $f(t_s)>0.2||(f(t_{j}-72),...,f(t_{j}))||_{\infty}$
\end{enumerate}
\item Determining an ending date based on discharge baseflow $Qb=Baseflow(Q)$:\\
$ed=\arg\min_{t_{e}}\sum_{t=t_{e}-1}^{t_{e}+47}|(Q-Qb)(t)|$ 
s.t. $t_{j} \leq t_e \leq sd+10\times24$
\end{enumerate}
\end{enumerate}
\textit{Remark.} If there exists $m+1$ $(m>0)$ consecutive events $(sd_{u},ed_{u}),...,(sd_{u+m},ed_{u+m})$
occurring ``nearly simultaneously'', that means all of these events
occur in no more than 10 days: $ed_{u+m}<sd_{u}+10\times24$, then we
merge these $m+1$ events into a single event $(sd_{u},ed_{u+m})$.}
\caption{Hydrograph segmentation algorithm}
\label{algo_seg}
\end{algorithm}

\section{Multi-objective optimization with spatially uniform control vectors}\label{MOO_appendix}

We look into multi-objective optimization for a global calibration of spatially uniform parameters, i.e. a low dimensional control $\overline{\boldsymbol{\theta}}$. The multi-objective calibration is simply
defined as the optimization problem:
\begin{equation}
\min_{\boldsymbol{\theta}\in\mathcal{O\subset\mathbb{R}}^{n}}(j_{1}(\boldsymbol{\theta}),...,j_{m}(\boldsymbol{\theta}))\label{prl:min}
\end{equation}
\noindent where $\boldsymbol{\theta}$ is the $n$-dimensional vector of model
parameters in the feasible space $\mathcal{\mathcal{O\subset\mathbb{R}}}^{n}$
and $j_{1},...,j_{m}$ are the $m$ single-objective functions to
be simultaneously minimized.
\subsection{Pareto front}\label{pareto_front_appendix}
\begin{figure}[ht!]
\begin{centering}
\includegraphics[scale=0.5]{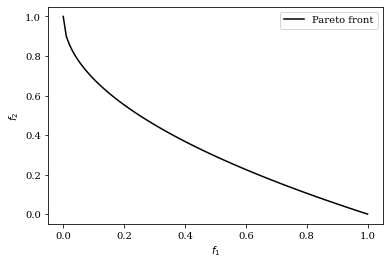}
\par\end{centering}
\caption{Illustration of Pareto front in terms of objective space.}
\label{pareto_front}
\end{figure}

In single-objective optimization, the Pareto optimal solution is unique
(in terms of objective space) but in multi-objective problem, it common
to have several solutions that cannot be defined which one is the
best. If the optimization problem is non-dominated, or non-inferior
(each objective function is its own entity, so no individual can be
better off without making at least one individual worse off), then
we call that Pareto optimality, or Pareto efficiency. A Pareto front
(in terms of parameter space) is a set of all Pareto efficient solutions
that need to be estimated. Let us consider two feasible solutions:
$\boldsymbol{\theta}_{1},\boldsymbol{\theta}_{2}\in\mathcal{O}$. Then, $\boldsymbol{\theta}_{1}$ is said
to Pareto dominate $\boldsymbol{\theta}_{2}$ if the following properties hold: 
\begin{enumerate}
\item $\forall i\in\{1,...,m\},j_{i}(\boldsymbol{\theta}_{1})\leq j_{i}(\boldsymbol{\theta}_{2})$;
\item $\exists i\in\{1,...,m\},j_{i}(\boldsymbol{\theta}_{1})<j_{i}(\boldsymbol{\theta}_{2})$.
\end{enumerate}
\noindent We call $\mathcal{P}$ the Pareto set representing all
of Pareto solutions. By definition, a Pareto solution
$\boldsymbol{\theta}^{*}\in\mathcal{P}$ of problem \ref{prl:min} must fill the two following conditions:
\begin{enumerate}
\item $\ensuremath{\nexists\boldsymbol{\theta}'\in}\mathcal{O}\setminus\mathcal{P},\exists i\in\{1,...,m\},j_{i}(\boldsymbol{\theta}')<j_{i}(\boldsymbol{\theta}^{*})$;
\item $\nexists\boldsymbol{\theta}''\in\mathcal{P},\boldsymbol{\theta}''$ dominates $\boldsymbol{\theta}^{*}$.
\end{enumerate}
The first statement indicates that there does not exist
other point in the feasible space that reduces at least one objective
function while keeping others unchanged, so the Pareto set is the
optimal set. The second says that, no other point exists in the Pareto
set that decreases one objective function without increasing another
one, so it is impossible to distinguish any solution as being better
than the other in the Pareto set. Fig. \ref{pareto_front} illustrates
this for a simple problem where we have 2-objective functions $j_{1}$,
$j_{2}$. The Pareto front (in terms of objective space) represents
all of non-dominated optimal solutions. It implies that, it is impossible
to move from any point in the feasible space and simultaneously decrease
the two objective functions without violating a constraint.
\subsection{Overview of GA and NSGA}\label{GA_overview}
GA is a ``heuristic algorithm'' (or search heuristic) in optimization, inspired by the Theory of Natural Evolution, whose selection operators include ``crossover'' and ``mutation''. Basically, the process of a GA consists of the following three phases:
\begin{enumerate}
    \item \textit{Population initialization.} The population is randomly initialized based on the problem range and constraint. The size of the population determines also the number of solutions, called ``pop-size''. 
    \item \textit{Parents selection (sorting).} A fitness function is defined to calculate the fitness score (also called Pareto ranking in multi-objective optimization) that determines how fit an individual is to the problem. Then, the fitness score decides the probability of selecting an individual as a parent to reproduce offspring population.
    \item \textit{Mating.} For each pair of parent to be mated, new offspring are created by exchanging the genes of parents among themselves (crossover operator). To maintain the diversity within the population and prevent premature convergence, some of the bits in the gene of certain new offspring can be flipped with a low random probability (mutation operator). Offspring are created until their pop-size is equal to the pop-size of previous generation. 
\end{enumerate}

\change{Suggested
by \citet{deb2002fast}, NSGA is a well known multi-objective genetic algorithm for solving multi/many-objective optimization problems, including fast and elitist approach \citep{deb2002fast}. Namely, a fast sorting algorithm helps optimizing the computational complexity (even with a large population size) arising from the non-dominated sorting procedure in every generation. Into the bargain, NSGA possesses a diversity preservation property, based on a sharing function method, that prevents the loss of good solutions involved in the mating process.
Recently, NSGA has also been implemented in the \textit{pymoo} Python library \citep{pymoo}, that is used in the present study thanks to the Python interface of our SMASH platform.}

\subsection{Selection of an optimal solution from Pareto front}\label{select_sol_append}
We aim to select an optimal solution that is acceptable for every objective within a constraint on principal objective function. Many strategies can be chosen to perform such a selection (e.g. based on the sensitivity ratio that is the ratio of the average variabilities of a certain non-inferior solution to the corresponding value of the objective function in the Pareto front \citep{wang2017multi}, or the Euclidean distance from the ideal solution \citep{wang2017application}). A simple additive weighting (SAW) method in \citet{wang2017application} can be used in our case by adding a normalization operator and assigned weightage for the objective functions.

Considering an objective matrix $(j_{ij})_{1\leq i\leq m, 1\leq j\leq n}$, where $m$ is the number of non-dominated solutions, $n$ is the number of objective functions. Then each row $i$ represents the $i^{th}$ solution set of the Pareto front and each column $j$ represents all non-inferior solutions of the $j^{th}$ objective function. Denote $c$ be the index of the classical objective function (for example $1-NSE$ or $1-KGE$), which is the most constrained function to find a unique optimal solution from Pareto front. This algorithm is detailed in the following three phases:
\begin{enumerate}
    \item Computing normalized objective matrix $(F_{ij})_{1\leq i\leq m, 1\leq j\leq n}$:
    $F_{ij}=\frac{f^{+}_j-f_{ij}}{f^{+}_j-f^{-}_j}$ \\
    where
    $f^{+}_j=\max_{1\leq i\leq m}f_{ij}$
    and
    $f^{-}_j=\min_{1\leq i\leq m}f_{ij}$
    \item Assigning weightage for normalized objective matrix $(G_{ij})_{1\leq i\leq m, 1\leq j\leq n}$:
    $G_{ij} = w_j \times F_{ij}$ \\
    where
        $w_j = \begin{cases}
            e^d \text{, if } j=c\\
            e - e^d \text{, otherwise.}
        \end{cases}
        $
        and
        $d = f^{+}_c-f^{-}_c$
    \item Finding optimal solution:
    $\boldsymbol{\theta} = (f_{k1}, ..., f_{kn})$
    where
    $k = \arg\max_{1\leq i\leq m} \left(\sum_{j=1}^n G_{ij}\right)$.
\end{enumerate}

\section{Calibration bounds}\label{BC_6param}

The parameter vector of SMASH model structure S6 is $\boldsymbol{\theta}(x)\equiv\left(c_{i}(x),c_{p}(x),c_{tr}(x),c_{r}(x),ml(x),c_{tl}(x)\right)^{T}$ and bound constrains used in optimization (Eq. \ref{inv_problem}) are set with values given in Table \ref{tab:BC-6-param}.

\begin{table}[ht!] \renewcommand*{\arraystretch}{1.25}
\begin{centering}
\caption{Boundary conditions of SMASH 6-parameters model.}
\label{tab:BC-6-param}
\scalebox{0.8}{
\begin{tabular}{c c c c c c c}
    {} & $c_i$ & $c_p$ & $c_{tr}$ & $c_{r}$ & $ml$ & $c_{tl}$\\
    \hline\hline
    Lower boundary & 1 & 1 & 1 & 1 & -20 & 1\\
    Upper boundary & 100 & 2000 & 1000 & 200 & 5 & 10000\\
\end{tabular}}
\par\end{centering}
\end{table}

\section*{Acknowledgments}
The authors greatly acknowledge SCHAPI-DGPR and METEO France for providing data used in this study. \oldrev{The first author was supported by funding from SCHAPI-DGPR and ANR grant ANR-21-CE04-0021-01 (MUFFINS project, "MUltiscale Flood Forecasting with INnovating
Solutions").} 

\printcredits

\bibliographystyle{elsarticle-harv}

\bibliography{references}

\begin{thebibliography}{86}
\expandafter\ifx\csname natexlab\endcsname\relax\def\natexlab#1{#1}\fi
\providecommand{\url}[1]{\texttt{#1}}
\providecommand{\href}[2]{#2}
\providecommand{\path}[1]{#1}
\providecommand{\DOIprefix}{doi:}
\providecommand{\ArXivprefix}{arXiv:}
\providecommand{\URLprefix}{URL: }
\providecommand{\Pubmedprefix}{pmid:}
\providecommand{\doi}[1]{\href{http://dx.doi.org/#1}{\path{#1}}}
\providecommand{\Pubmed}[1]{\href{pmid:#1}{\path{#1}}}
\providecommand{\bibinfo}[2]{#2}
\ifx\xfnm\relax \def\xfnm[#1]{\unskip,\space#1}\fi
\bibitem[{Astagneau et~al.(2021)Astagneau, Bourgin, Andréassian and
  Perrin}]{astagneau_2021}
\bibinfo{author}{Astagneau, P.C.}, \bibinfo{author}{Bourgin, F.},
  \bibinfo{author}{Andréassian, V.}, \bibinfo{author}{Perrin, C.},
  \bibinfo{year}{2021}.
\newblock \bibinfo{title}{When does a parsimonious model fail to simulate
  floods? learning from the seasonality of model bias}.
\newblock \bibinfo{journal}{Hydrological Sciences Journal}
  \bibinfo{volume}{66}, \bibinfo{pages}{1288--1305}.
\newblock \URLprefix \url{https://doi.org/10.1080/02626667.2021.1923720},
  \DOIprefix\doi{10.1080/02626667.2021.1923720},
  \href{http://arxiv.org/abs/https://doi.org/10.1080/02626667.2021.1923720}{{\tt
  arXiv:https://doi.org/10.1080/02626667.2021.1923720}}.
\bibitem[{Azzini et~al.(2021)Azzini, Mara and Rosati}]{azzini2021comparison}
\bibinfo{author}{Azzini, I.}, \bibinfo{author}{Mara, T.A.},
  \bibinfo{author}{Rosati, R.}, \bibinfo{year}{2021}.
\newblock \bibinfo{title}{Comparison of two sets of monte carlo estimators of
  sobol’indices}.
\newblock \bibinfo{journal}{Environmental Modelling \& Software}
  \bibinfo{volume}{144}, \bibinfo{pages}{105167}.
\bibitem[{Beck(1987)}]{beck1987water}
\bibinfo{author}{Beck, M.B.}, \bibinfo{year}{1987}.
\newblock \bibinfo{title}{Water quality modeling: a review of the analysis of
  uncertainty}.
\newblock \bibinfo{journal}{Water resources research} \bibinfo{volume}{23},
  \bibinfo{pages}{1393--1442}.
\bibitem[{Beven(1993)}]{beven1993prophecy}
\bibinfo{author}{Beven, K.}, \bibinfo{year}{1993}.
\newblock \bibinfo{title}{Prophecy, reality and uncertainty in distributed
  hydrological modelling}.
\newblock \bibinfo{journal}{Advances in water resources} \bibinfo{volume}{16},
  \bibinfo{pages}{41--51}.
\bibitem[{{Blank} and {Deb}(2020)}]{pymoo}
\bibinfo{author}{{Blank}, J.}, \bibinfo{author}{{Deb}, K.},
  \bibinfo{year}{2020}.
\newblock \bibinfo{title}{pymoo: Multi-objective optimization in python}.
\newblock \bibinfo{journal}{IEEE Access} \bibinfo{volume}{8},
  \bibinfo{pages}{89497--89509}.
\bibitem[{Bout and Jetten(2018)}]{bout2018validity}
\bibinfo{author}{Bout, V.}, \bibinfo{author}{Jetten, V.}, \bibinfo{year}{2018}.
\newblock \bibinfo{title}{The validity of flow approximations when simulating
  catchment-integrated flash floods}.
\newblock \bibinfo{journal}{Journal of hydrology} \bibinfo{volume}{556},
  \bibinfo{pages}{674--688}.
\bibitem[{Brigode et~al.(2020)Brigode, G{\'e}not, Lobligeois and
  Delaigue}]{brigode2020data}
\bibinfo{author}{Brigode, P.}, \bibinfo{author}{G{\'e}not, B.},
  \bibinfo{author}{Lobligeois, F.}, \bibinfo{author}{Delaigue, O.},
  \bibinfo{year}{2020}.
\newblock \bibinfo{title}{Summary sheets of watershed-scale hydroclimatic
  observed data for france}.
\newblock \URLprefix \url{https://doi.org/10.15454/UV01P1},
  \DOIprefix\doi{10.15454/UV01P1}.
\bibitem[{Campolongo et~al.(2011)Campolongo, Saltelli and
  Cariboni}]{campolongo2011screening}
\bibinfo{author}{Campolongo, F.}, \bibinfo{author}{Saltelli, A.},
  \bibinfo{author}{Cariboni, J.}, \bibinfo{year}{2011}.
\newblock \bibinfo{title}{From screening to quantitative sensitivity analysis.
  a unified approach}.
\newblock \bibinfo{journal}{Computer physics communications}
  \bibinfo{volume}{182}, \bibinfo{pages}{978--988}.
\bibitem[{Champeaux et~al.(2009)Champeaux, Dupuy, Laurantin, Soulan, Tabary and
  Soubeyroux}]{champeaux2009mesures}
\bibinfo{author}{Champeaux, J.L.}, \bibinfo{author}{Dupuy, P.},
  \bibinfo{author}{Laurantin, O.}, \bibinfo{author}{Soulan, I.},
  \bibinfo{author}{Tabary, P.}, \bibinfo{author}{Soubeyroux, J.M.},
  \bibinfo{year}{2009}.
\newblock \bibinfo{title}{Les mesures de pr{\'e}cipitations et l'estimation des
  lames d'eau {\`a} m{\'e}t{\'e}o-france: {\'e}tat de l'art et perspectives}.
\newblock \bibinfo{journal}{La Houille Blanche} , \bibinfo{pages}{28--34}.
\bibitem[{Chibeles-Martins et~al.(2016)Chibeles-Martins, Pinto-Varela,
  Barbosa-P{\'o}voa and Novais}]{chibeles2016multi}
\bibinfo{author}{Chibeles-Martins, N.}, \bibinfo{author}{Pinto-Varela, T.},
  \bibinfo{author}{Barbosa-P{\'o}voa, A.P.}, \bibinfo{author}{Novais, A.Q.},
  \bibinfo{year}{2016}.
\newblock \bibinfo{title}{A multi-objective meta-heuristic approach for the
  design and planning of green supply chains-mbsa}.
\newblock \bibinfo{journal}{Expert Systems with Applications}
  \bibinfo{volume}{47}, \bibinfo{pages}{71--84}.
\bibitem[{Colleoni et~al.(2022)Colleoni, Garambois, Javelle, Jay-Allemand and
  Arnaud}]{egusphere-2022-506}
\bibinfo{author}{Colleoni, F.}, \bibinfo{author}{Garambois, P.A.},
  \bibinfo{author}{Javelle, P.}, \bibinfo{author}{Jay-Allemand, M.},
  \bibinfo{author}{Arnaud, P.}, \bibinfo{year}{2022}.
\newblock \bibinfo{title}{Adjoint-based spatially distributed calibration of a
  grid gr-based parsimonious hydrological model over 312 french catchments with
  smash platform}.
\newblock \bibinfo{journal}{EGUsphere} \bibinfo{volume}{2022},
  \bibinfo{pages}{1--37}.
\newblock \URLprefix
  \url{https://egusphere.copernicus.org/preprints/2022/egusphere-2022-506/},
  \DOIprefix\doi{10.5194/egusphere-2022-506}.
\bibitem[{De~Lavenne et~al.(2019)De~Lavenne, Andr{\'e}assian, Thirel, Ramos and
  Perrin}]{de2019regularization}
\bibinfo{author}{De~Lavenne, A.}, \bibinfo{author}{Andr{\'e}assian, V.},
  \bibinfo{author}{Thirel, G.}, \bibinfo{author}{Ramos, M.H.},
  \bibinfo{author}{Perrin, C.}, \bibinfo{year}{2019}.
\newblock \bibinfo{title}{A regularization approach to improve the sequential
  calibration of a semidistributed hydrological model}.
\newblock \bibinfo{journal}{Water Resources Research} \bibinfo{volume}{55},
  \bibinfo{pages}{8821--8839}.
\bibitem[{Deb et~al.(2002)Deb, Pratap, Agarwal and Meyarivan}]{deb2002fast}
\bibinfo{author}{Deb, K.}, \bibinfo{author}{Pratap, A.},
  \bibinfo{author}{Agarwal, S.}, \bibinfo{author}{Meyarivan, T.},
  \bibinfo{year}{2002}.
\newblock \bibinfo{title}{A fast and elitist multiobjective genetic algorithm:
  Nsga-ii}.
\newblock \bibinfo{journal}{IEEE transactions on evolutionary computation}
  \bibinfo{volume}{6}, \bibinfo{pages}{182--197}.
\bibitem[{Delaigue et~al.(2020)Delaigue, G{\'e}not, Lebecherel, Brigode and
  Bourgin}]{delaigue2020data}
\bibinfo{author}{Delaigue, O.}, \bibinfo{author}{G{\'e}not, B.},
  \bibinfo{author}{Lebecherel, L.}, \bibinfo{author}{Brigode, P.},
  \bibinfo{author}{Bourgin, P.Y.}, \bibinfo{year}{2020}.
\newblock \bibinfo{title}{Database of watershed-scale hydroclimatic
  observations in france}.
\newblock \URLprefix \url{https://webgr.inrae.fr/base-de-donnees}.
\bibitem[{Douinot et~al.(2018)Douinot, Roux, Garambois and
  Dartus}]{douinot2018using}
\bibinfo{author}{Douinot, A.}, \bibinfo{author}{Roux, H.},
  \bibinfo{author}{Garambois, P.A.}, \bibinfo{author}{Dartus, D.},
  \bibinfo{year}{2018}.
\newblock \bibinfo{title}{Using a multi-hypothesis framework to improve the
  understanding of flow dynamics during flash floods}.
\newblock \bibinfo{journal}{Hydrology and Earth System Sciences}
  \bibinfo{volume}{22}, \bibinfo{pages}{5317--5340}.
\bibitem[{Duarte and Watanabe(2021)}]{marcos_duarte_2021_4599319}
\bibinfo{author}{Duarte, M.}, \bibinfo{author}{Watanabe, R.N.},
  \bibinfo{year}{2021}.
\newblock \bibinfo{title}{{Notes on Scientific Computing for Biomechanics and
  Motor Control}}.
\newblock \URLprefix \url{https://doi.org/10.5281/zenodo.4599319},
  \DOIprefix\doi{10.5281/zenodo.4599319}.
\bibitem[{Edijatno(1991)}]{edijatno91}
\bibinfo{author}{Edijatno}, \bibinfo{year}{1991}.
\newblock \bibinfo{title}{Mise au point d'un modele elementaire pluie-debit au
  pas de temps journalier}.
\newblock Ph.D. thesis. Universite Louis Pasteur, ENGEES, Cemagref Antony.
\bibitem[{El-Ghandour and Elbeltagi(2014)}]{el2014optimal}
\bibinfo{author}{El-Ghandour, H.A.}, \bibinfo{author}{Elbeltagi, E.},
  \bibinfo{year}{2014}.
\newblock \bibinfo{title}{Optimal groundwater management using multiobjective
  particle swarm with a new evolution strategy}.
\newblock \bibinfo{journal}{Journal of Hydrologic Engineering}
  \bibinfo{volume}{19}, \bibinfo{pages}{1141--1149}.
\bibitem[{Emmanuel et~al.(2015)Emmanuel, Andrieu, Leblois, Janey and
  Payrastre}]{emmanuel2015influence}
\bibinfo{author}{Emmanuel, I.}, \bibinfo{author}{Andrieu, H.},
  \bibinfo{author}{Leblois, E.}, \bibinfo{author}{Janey, N.},
  \bibinfo{author}{Payrastre, O.}, \bibinfo{year}{2015}.
\newblock \bibinfo{title}{Influence of rainfall spatial variability on
  rainfall--runoff modelling: benefit of a simulation approach?}
\newblock \bibinfo{journal}{Journal of hydrology} \bibinfo{volume}{531},
  \bibinfo{pages}{337--348}.
\bibitem[{Garambois et~al.(2014)Garambois, Larnier, Roux, Labat and
  Dartus}]{garambois2014analysis}
\bibinfo{author}{Garambois, P.A.}, \bibinfo{author}{Larnier, K.},
  \bibinfo{author}{Roux, H.}, \bibinfo{author}{Labat, D.},
  \bibinfo{author}{Dartus, D.}, \bibinfo{year}{2014}.
\newblock \bibinfo{title}{Analysis of flash flood-triggering rainfall for a
  process-oriented hydrological model}.
\newblock \bibinfo{journal}{Atmospheric research} \bibinfo{volume}{137},
  \bibinfo{pages}{14--24}.
\bibitem[{Garambois et~al.(2013)Garambois, Roux, Larnier, Castaings and
  Dartus}]{garambois2013sobol}
\bibinfo{author}{Garambois, P.A.}, \bibinfo{author}{Roux, H.},
  \bibinfo{author}{Larnier, K.}, \bibinfo{author}{Castaings, W.},
  \bibinfo{author}{Dartus, D.}, \bibinfo{year}{2013}.
\newblock \bibinfo{title}{Characterization of process-oriented hydrologic model
  behavior with temporal sensitivity analysis for flash floods in mediterranean
  catchments}.
\newblock \bibinfo{journal}{Hydrology and Earth System Sciences} .
\bibitem[{Garambois et~al.(2015)Garambois, Roux, Larnier, Labat and
  Dartus}]{garambois2015characterization}
\bibinfo{author}{Garambois, P.A.}, \bibinfo{author}{Roux, H.},
  \bibinfo{author}{Larnier, K.}, \bibinfo{author}{Labat, D.},
  \bibinfo{author}{Dartus, D.}, \bibinfo{year}{2015}.
\newblock \bibinfo{title}{Characterization of catchment behaviour and rainfall
  selection for flash flood hydrological model calibration: catchments of the
  eastern pyrenees}.
\newblock \bibinfo{journal}{Hydrological sciences journal}
  \bibinfo{volume}{60}, \bibinfo{pages}{424--447}.
\bibitem[{Guo et~al.(2014)Guo, Zhou, Lu, Zou, Zhang and Bi}]{guo2014multi}
\bibinfo{author}{Guo, J.}, \bibinfo{author}{Zhou, J.}, \bibinfo{author}{Lu,
  J.}, \bibinfo{author}{Zou, Q.}, \bibinfo{author}{Zhang, H.},
  \bibinfo{author}{Bi, S.}, \bibinfo{year}{2014}.
\newblock \bibinfo{title}{Multi-objective optimization of empirical
  hydrological model for streamflow prediction}.
\newblock \bibinfo{journal}{Journal of Hydrology} \bibinfo{volume}{511},
  \bibinfo{pages}{242--253}.
\bibitem[{Gupta et~al.(2006)Gupta, Beven and Wagener}]{gupta2006model}
\bibinfo{author}{Gupta, H.V.}, \bibinfo{author}{Beven, K.J.},
  \bibinfo{author}{Wagener, T.}, \bibinfo{year}{2006}.
\newblock \bibinfo{title}{Model calibration and uncertainty estimation}.
\newblock \bibinfo{journal}{Encyclopedia of hydrological sciences} .
\bibitem[{Gupta et~al.(2009)Gupta, Kling, Yilmaz and
  Martinez}]{gupta2009decomposition}
\bibinfo{author}{Gupta, H.V.}, \bibinfo{author}{Kling, H.},
  \bibinfo{author}{Yilmaz, K.K.}, \bibinfo{author}{Martinez, G.F.},
  \bibinfo{year}{2009}.
\newblock \bibinfo{title}{Decomposition of the mean squared error and nse
  performance criteria: Implications for improving hydrological modelling}.
\newblock \bibinfo{journal}{Journal of hydrology} \bibinfo{volume}{377},
  \bibinfo{pages}{80--91}.
\bibitem[{Gupta and Razavi(2018)}]{gupta2018revisiting}
\bibinfo{author}{Gupta, H.V.}, \bibinfo{author}{Razavi, S.},
  \bibinfo{year}{2018}.
\newblock \bibinfo{title}{Revisiting the basis of sensitivity analysis for
  dynamical earth system models}.
\newblock \bibinfo{journal}{Water Resources Research} \bibinfo{volume}{54},
  \bibinfo{pages}{8692--8717}.
\bibitem[{Gupta et~al.(1998)Gupta, Sorooshian and Yapo}]{gupta1998toward}
\bibinfo{author}{Gupta, H.V.}, \bibinfo{author}{Sorooshian, S.},
  \bibinfo{author}{Yapo, P.O.}, \bibinfo{year}{1998}.
\newblock \bibinfo{title}{Toward improved calibration of hydrologic models:
  Multiple and noncommensurable measures of information}.
\newblock \bibinfo{journal}{Water Resources Research} \bibinfo{volume}{34},
  \bibinfo{pages}{751--763}.
\bibitem[{Hascoet and Pascual(2013)}]{tapenade}
\bibinfo{author}{Hascoet, L.}, \bibinfo{author}{Pascual, V.},
  \bibinfo{year}{2013}.
\newblock \bibinfo{title}{The tapenade automatic differentiation tool:
  Principles, model, and specification}.
\newblock \bibinfo{journal}{ACM Trans. Math. Softw.} \bibinfo{volume}{39}.
\newblock \URLprefix \url{https://doi.org/10.1145/2450153.2450158},
  \DOIprefix\doi{10.1145/2450153.2450158}.
\bibitem[{Herman and Usher(2017)}]{Herman2017}
\bibinfo{author}{Herman, J.}, \bibinfo{author}{Usher, W.},
  \bibinfo{year}{2017}.
\newblock \bibinfo{title}{{SALib}: An open-source python library for
  sensitivity analysis}.
\newblock \bibinfo{journal}{The Journal of Open Source Software}
  \bibinfo{volume}{2}.
\newblock \URLprefix \url{https://doi.org/10.21105/joss.00097},
  \DOIprefix\doi{10.21105/joss.00097}.
\bibitem[{Horner(2020)}]{horner_PhD_2020}
\bibinfo{author}{Horner, I.}, \bibinfo{year}{2020}.
\newblock \bibinfo{title}{Design and Evaluation of hydrological signatures for
  the diagnosis and improvement of a process-based distributed hydrological
  model}.
\newblock Ph.D. thesis. Université Grenoble Alpes.
\bibitem[{Iooss and Lema{\^\i}tre(2015)}]{iooss2015review}
\bibinfo{author}{Iooss, B.}, \bibinfo{author}{Lema{\^\i}tre, P.},
  \bibinfo{year}{2015}.
\newblock \bibinfo{title}{A review on global sensitivity analysis methods}.
\newblock \bibinfo{journal}{Uncertainty management in simulation-optimization
  of complex systems: algorithms and applications} , \bibinfo{pages}{101--122}.
\bibitem[{Iwanaga et~al.(2022)Iwanaga, Usher and Herman}]{Iwanaga2022}
\bibinfo{author}{Iwanaga, T.}, \bibinfo{author}{Usher, W.},
  \bibinfo{author}{Herman, J.}, \bibinfo{year}{2022}.
\newblock \bibinfo{title}{Toward {SALib} 2.0: {Advancing} the accessibility and
  interpretability of global sensitivity analyses}.
\newblock \bibinfo{journal}{Socio-Environmental Systems Modelling}
  \bibinfo{volume}{4}, \bibinfo{pages}{18155}.
\newblock \URLprefix \url{https://sesmo.org/article/view/18155},
  \DOIprefix\doi{10.18174/sesmo.18155}.
\bibitem[{Jay-Allemand(2020)}]{jay2020estimation}
\bibinfo{author}{Jay-Allemand, M.}, \bibinfo{year}{2020}.
\newblock \bibinfo{title}{Estimation variationnelle des param{\`e}tres d'un
  mod{\`e}le hydrologique distribu{\'e}}.
\newblock Ph.D. thesis. Aix-Marseille.
\bibitem[{Jay-Allemand et~al.(2022a)Jay-Allemand, Colleoni, Garambois, Javelle
  and Julie}]{JayCoGa_2022_wrappingSMASH}
\bibinfo{author}{Jay-Allemand, M.}, \bibinfo{author}{Colleoni, F.},
  \bibinfo{author}{Garambois, P.A.}, \bibinfo{author}{Javelle, P.},
  \bibinfo{author}{Julie, D.}, \bibinfo{year}{2022}a.
\newblock \bibinfo{title}{Smash - spatially distributed modelling and
  assimilation for hydrology: Python wrapping towards enhanced
  research-to-operations transfer}.
\newblock \bibinfo{journal}{IAHS} \URLprefix
  \url{https://hal.archives-ouvertes.fr/hal-03683657}.
\bibitem[{Jay-Allemand et~al.(2022b)Jay-Allemand, Demargne, Garambois, Javelle,
  Gejadze, Colleoni, Organde, Arnaud and Fouchier}]{jay2022spatially}
\bibinfo{author}{Jay-Allemand, M.}, \bibinfo{author}{Demargne, J.},
  \bibinfo{author}{Garambois, P.A.}, \bibinfo{author}{Javelle, P.},
  \bibinfo{author}{Gejadze, I.}, \bibinfo{author}{Colleoni, F.},
  \bibinfo{author}{Organde, D.}, \bibinfo{author}{Arnaud, P.},
  \bibinfo{author}{Fouchier, C.}, \bibinfo{year}{2022}b.
\newblock \bibinfo{title}{Spatially distributed calibration of a hydrological
  model with variational optimization constrained by physiographic maps for
  flash flood forecasting in France}.
\newblock \bibinfo{type}{Technical Report}. Copernicus Meetings.
\bibitem[{Jay-Allemand et~al.(2020)Jay-Allemand, Javelle, Gejadze, Arnaud,
  Malaterre, Fine and Organde}]{jay2020potential}
\bibinfo{author}{Jay-Allemand, M.}, \bibinfo{author}{Javelle, P.},
  \bibinfo{author}{Gejadze, I.}, \bibinfo{author}{Arnaud, P.},
  \bibinfo{author}{Malaterre, P.O.}, \bibinfo{author}{Fine, J.A.},
  \bibinfo{author}{Organde, D.}, \bibinfo{year}{2020}.
\newblock \bibinfo{title}{On the potential of variational calibration for a
  fully distributed hydrological model: application on a mediterranean
  catchment}.
\newblock \bibinfo{journal}{Hydrology and Earth System Sciences}
  \bibinfo{volume}{24}, \bibinfo{pages}{5519--5538}.
\bibitem[{Kavetski et~al.(2018)Kavetski, Fenicia, Reichert and
  Albert}]{Kavetski_2017WR}
\bibinfo{author}{Kavetski, D.}, \bibinfo{author}{Fenicia, F.},
  \bibinfo{author}{Reichert, P.}, \bibinfo{author}{Albert, C.},
  \bibinfo{year}{2018}.
\newblock \bibinfo{title}{Signature-domain calibration of hydrological models
  using approximate bayesian computation: Theory and comparison to existing
  applications}.
\newblock \bibinfo{journal}{Water Resources Research} \bibinfo{volume}{54},
  \bibinfo{pages}{4059--4083}.
\newblock \URLprefix
  \url{https://agupubs.onlinelibrary.wiley.com/doi/abs/10.1002/2017WR020528},
  \DOIprefix\doi{https://doi.org/10.1002/2017WR020528},
  \href{http://arxiv.org/abs/https://agupubs.onlinelibrary.wiley.com/doi/pdf/10.1002/2017WR020528}{{\tt
  arXiv:https://agupubs.onlinelibrary.wiley.com/doi/pdf/10.1002/2017WR020528}}.
\bibitem[{Khorram et~al.(2014)Khorram, Khaledian and
  Khaledyan}]{khorram2014numerical}
\bibinfo{author}{Khorram, E.}, \bibinfo{author}{Khaledian, K.},
  \bibinfo{author}{Khaledyan, M.}, \bibinfo{year}{2014}.
\newblock \bibinfo{title}{A numerical method for constructing the pareto front
  of multi-objective optimization problems}.
\newblock \bibinfo{journal}{Journal of Computational and Applied Mathematics}
  \bibinfo{volume}{261}, \bibinfo{pages}{158--171}.
\bibitem[{Kirstetter et~al.(2021)Kirstetter, Delestre, Lagr{\'e}e, Popinet and
  Josserand}]{kirstetter2021b}
\bibinfo{author}{Kirstetter, G.}, \bibinfo{author}{Delestre, O.},
  \bibinfo{author}{Lagr{\'e}e, P.Y.}, \bibinfo{author}{Popinet, S.},
  \bibinfo{author}{Josserand, C.}, \bibinfo{year}{2021}.
\newblock \bibinfo{title}{B-flood 1.0: an open-source saint-venant model for
  flash-flood simulation using adaptive refinement}.
\newblock \bibinfo{journal}{Geoscientific Model Development}
  \bibinfo{volume}{14}, \bibinfo{pages}{7117--7132}.
\bibitem[{Lamboni et~al.(2013)Lamboni, Iooss, Popelin and
  Gamboa}]{lamboni2013derivative}
\bibinfo{author}{Lamboni, M.}, \bibinfo{author}{Iooss, B.},
  \bibinfo{author}{Popelin, A.L.}, \bibinfo{author}{Gamboa, F.},
  \bibinfo{year}{2013}.
\newblock \bibinfo{title}{Derivative-based global sensitivity measures: General
  links with sobol’indices and numerical tests}.
\newblock \bibinfo{journal}{Mathematics and Computers in Simulation}
  \bibinfo{volume}{87}, \bibinfo{pages}{45--54}.
\bibitem[{Le~Mesnil(2021)}]{lemesnil:tel-03578569}
\bibinfo{author}{Le~Mesnil, M.}, \bibinfo{year}{2021}.
\newblock \bibinfo{title}{{Signatures Hydrologiques des Bassins Karstiques}}.
\newblock \bibinfo{type}{Theses}. {Montpellier SupAgro}.
\newblock \URLprefix \url{https://tel.archives-ouvertes.fr/tel-03578569}.
\bibitem[{Li et~al.(2022)Li, Huang, Wang, Razavi and Zhang}]{li2022development}
\bibinfo{author}{Li, K.}, \bibinfo{author}{Huang, G.}, \bibinfo{author}{Wang,
  S.}, \bibinfo{author}{Razavi, S.}, \bibinfo{author}{Zhang, X.},
  \bibinfo{year}{2022}.
\newblock \bibinfo{title}{Development of a joint probabilistic rainfall-runoff
  model for high-to-extreme flow projections under changing climatic
  conditions}.
\newblock \bibinfo{journal}{Water Resources Research} \bibinfo{volume}{58},
  \bibinfo{pages}{e2021WR031557}.
\bibitem[{Li et~al.(2021)Li, Chen, Gao, Luo, Gourley, Kirstetter, Yang, Kolar,
  McGovern, Wen et~al.}]{li2021crest}
\bibinfo{author}{Li, Z.}, \bibinfo{author}{Chen, M.}, \bibinfo{author}{Gao,
  S.}, \bibinfo{author}{Luo, X.}, \bibinfo{author}{Gourley, J.J.},
  \bibinfo{author}{Kirstetter, P.}, \bibinfo{author}{Yang, T.},
  \bibinfo{author}{Kolar, R.}, \bibinfo{author}{McGovern, A.},
  \bibinfo{author}{Wen, Y.}, et~al., \bibinfo{year}{2021}.
\newblock \bibinfo{title}{Crest-imap v1. 0: A fully coupled
  hydrologic-hydraulic modeling framework dedicated to flood inundation mapping
  and prediction}.
\newblock \bibinfo{journal}{Environmental Modelling \& Software}
  \bibinfo{volume}{141}, \bibinfo{pages}{105051}.
\bibitem[{Lyne and Hollick(1979)}]{lyne1979stochastic}
\bibinfo{author}{Lyne, V.}, \bibinfo{author}{Hollick, M.},
  \bibinfo{year}{1979}.
\newblock \bibinfo{title}{Stochastic time-variable rainfall-runoff modelling},
  in: \bibinfo{booktitle}{Institute of Engineers Australia National
  Conference}, \bibinfo{organization}{Institute of Engineers Australia Barton,
  Australia}. pp. \bibinfo{pages}{89--93}.
\bibitem[{McMillan(2021)}]{mcmillan2021review}
\bibinfo{author}{McMillan, H.K.}, \bibinfo{year}{2021}.
\newblock \bibinfo{title}{A review of hydrologic signatures and their
  applications}.
\newblock \bibinfo{journal}{Wiley Interdisciplinary Reviews: Water}
  \bibinfo{volume}{8}, \bibinfo{pages}{e1499}.
\bibitem[{Mizukami et~al.(2017)Mizukami, Clark, Newman, Wood, Gutmann, Nijssen,
  Rakovec and Samaniego}]{mizukami2017towards}
\bibinfo{author}{Mizukami, N.}, \bibinfo{author}{Clark, M.P.},
  \bibinfo{author}{Newman, A.J.}, \bibinfo{author}{Wood, A.W.},
  \bibinfo{author}{Gutmann, E.D.}, \bibinfo{author}{Nijssen, B.},
  \bibinfo{author}{Rakovec, O.}, \bibinfo{author}{Samaniego, L.},
  \bibinfo{year}{2017}.
\newblock \bibinfo{title}{Towards seamless large-domain parameter estimation
  for hydrologic models}.
\newblock \bibinfo{journal}{Water Resources Research} \bibinfo{volume}{53},
  \bibinfo{pages}{8020--8040}.
\bibitem[{Mizukami et~al.(2019)Mizukami, Rakovec, Newman, Clark, Wood, Gupta
  and Kumar}]{mizukami2019choice}
\bibinfo{author}{Mizukami, N.}, \bibinfo{author}{Rakovec, O.},
  \bibinfo{author}{Newman, A.J.}, \bibinfo{author}{Clark, M.P.},
  \bibinfo{author}{Wood, A.W.}, \bibinfo{author}{Gupta, H.V.},
  \bibinfo{author}{Kumar, R.}, \bibinfo{year}{2019}.
\newblock \bibinfo{title}{On the choice of calibration metrics for
  “high-flow” estimation using hydrologic models}.
\newblock \bibinfo{journal}{Hydrology and Earth System Sciences}
  \bibinfo{volume}{23}, \bibinfo{pages}{2601--2614}.
\bibitem[{Monnier et~al.(2016)Monnier, Couderc, Dartus, Larnier, Madec and
  Vila}]{monnier2016inverse}
\bibinfo{author}{Monnier, J.}, \bibinfo{author}{Couderc, F.},
  \bibinfo{author}{Dartus, D.}, \bibinfo{author}{Larnier, K.},
  \bibinfo{author}{Madec, R.}, \bibinfo{author}{Vila, J.P.},
  \bibinfo{year}{2016}.
\newblock \bibinfo{title}{Inverse algorithms for 2d shallow water equations in
  presence of wet dry fronts: Application to flood plain dynamics}.
\newblock \bibinfo{journal}{Advances in Water Resources} \bibinfo{volume}{97},
  \bibinfo{pages}{11--24}.
\bibitem[{Mostafaie et~al.(2018)Mostafaie, Forootan, Safari and
  Schumacher}]{mostafaie2018comparing}
\bibinfo{author}{Mostafaie, A.}, \bibinfo{author}{Forootan, E.},
  \bibinfo{author}{Safari, A.}, \bibinfo{author}{Schumacher, M.},
  \bibinfo{year}{2018}.
\newblock \bibinfo{title}{Comparing multi-objective optimization techniques to
  calibrate a conceptual hydrological model using in situ runoff and daily
  grace data}.
\newblock \bibinfo{journal}{Computational Geosciences} \bibinfo{volume}{22},
  \bibinfo{pages}{789--814}.
\bibitem[{Nash and Sutcliffe(1970)}]{nash1970river}
\bibinfo{author}{Nash, J.E.}, \bibinfo{author}{Sutcliffe, J.V.},
  \bibinfo{year}{1970}.
\newblock \bibinfo{title}{River flow forecasting through conceptual models part
  i—a discussion of principles}.
\newblock \bibinfo{journal}{Journal of hydrology} \bibinfo{volume}{10},
  \bibinfo{pages}{282--290}.
\bibitem[{Nathan and McMahon(1990)}]{nathan1990evaluation}
\bibinfo{author}{Nathan, R.J.}, \bibinfo{author}{McMahon, T.A.},
  \bibinfo{year}{1990}.
\newblock \bibinfo{title}{Evaluation of automated techniques for base flow and
  recession analyses}.
\newblock \bibinfo{journal}{Water resources research} \bibinfo{volume}{26},
  \bibinfo{pages}{1465--1473}.
\bibitem[{Oliveira et~al.(2021)Oliveira, Fleischmann and
  Paiva}]{oliveira2021contribution}
\bibinfo{author}{Oliveira, A.M.}, \bibinfo{author}{Fleischmann, A.},
  \bibinfo{author}{Paiva, R.}, \bibinfo{year}{2021}.
\newblock \bibinfo{title}{On the contribution of remote sensing-based
  calibration to model hydrological and hydraulic processes in tropical
  regions}.
\newblock \bibinfo{journal}{Journal of Hydrology} \bibinfo{volume}{597},
  \bibinfo{pages}{126184}.
\bibitem[{Oudin et~al.(2005)Oudin, Hervieu, Michel, Perrin, Andr{\'e}assian,
  Anctil and Loumagne}]{oudin2005potential}
\bibinfo{author}{Oudin, L.}, \bibinfo{author}{Hervieu, F.},
  \bibinfo{author}{Michel, C.}, \bibinfo{author}{Perrin, C.},
  \bibinfo{author}{Andr{\'e}assian, V.}, \bibinfo{author}{Anctil, F.},
  \bibinfo{author}{Loumagne, C.}, \bibinfo{year}{2005}.
\newblock \bibinfo{title}{Which potential evapotranspiration input for a lumped
  rainfall--runoff model?: Part 2 towards a simple and efficient potential
  evapotranspiration model for rainfall--runoff modelling}.
\newblock \bibinfo{journal}{Journal of hydrology} \bibinfo{volume}{303},
  \bibinfo{pages}{290--306}.
\bibitem[{Paiva et~al.(2011)Paiva, Collischonn and Tucci}]{paiva2011large}
\bibinfo{author}{Paiva, R.C.}, \bibinfo{author}{Collischonn, W.},
  \bibinfo{author}{Tucci, C.E.}, \bibinfo{year}{2011}.
\newblock \bibinfo{title}{Large scale hydrologic and hydrodynamic modeling
  using limited data and a gis based approach}.
\newblock \bibinfo{journal}{Journal of Hydrology} \bibinfo{volume}{406},
  \bibinfo{pages}{170--181}.
\bibitem[{Pelletier and Andr{\'e}assian(2020)}]{pelletier2020hydrograph}
\bibinfo{author}{Pelletier, A.}, \bibinfo{author}{Andr{\'e}assian, V.},
  \bibinfo{year}{2020}.
\newblock \bibinfo{title}{Hydrograph separation: an impartial parametrisation
  for an imperfect method}.
\newblock \bibinfo{journal}{Hydrology and earth system sciences}
  \bibinfo{volume}{24}, \bibinfo{pages}{1171--1187}.
\bibitem[{Pujol et~al.(2020)Pujol, Garambois, Finaud-Guyot, Monnier, Larnier,
  Mose, Biancamaria, Yesou, Moreira, Paris et~al.}]{pujol2020estimation}
\bibinfo{author}{Pujol, L.}, \bibinfo{author}{Garambois, P.A.},
  \bibinfo{author}{Finaud-Guyot, P.}, \bibinfo{author}{Monnier, J.},
  \bibinfo{author}{Larnier, K.}, \bibinfo{author}{Mose, R.},
  \bibinfo{author}{Biancamaria, S.}, \bibinfo{author}{Yesou, H.},
  \bibinfo{author}{Moreira, D.}, \bibinfo{author}{Paris, A.}, et~al.,
  \bibinfo{year}{2020}.
\newblock \bibinfo{title}{Estimation of multiple inflows and effective channel
  by assimilation of multi-satellite hydraulic signatures: The ungauged
  anabranching negro river}.
\newblock \bibinfo{journal}{Journal of Hydrology} \bibinfo{volume}{591},
  \bibinfo{pages}{125331}.
\bibitem[{Pujol et~al.(2022)Pujol, Garambois and Monnier}]{pujol2022multi}
\bibinfo{author}{Pujol, L.}, \bibinfo{author}{Garambois, P.A.},
  \bibinfo{author}{Monnier, J.}, \bibinfo{year}{2022}.
\newblock \bibinfo{title}{Multi-dimensional hydrological--hydraulic model with
  variational data assimilation for river networks and floodplains}.
\newblock \bibinfo{journal}{Geoscientific Model Development}
  \bibinfo{volume}{15}, \bibinfo{pages}{6085--6113}.
\bibitem[{Puy et~al.(2022)Puy, Becker, Piano and
  Saltelli}]{puy2022comprehensive}
\bibinfo{author}{Puy, A.}, \bibinfo{author}{Becker, W.},
  \bibinfo{author}{Piano, S.L.}, \bibinfo{author}{Saltelli, A.},
  \bibinfo{year}{2022}.
\newblock \bibinfo{title}{A comprehensive comparison of total-order estimators
  for global sensitivity analysis}.
\newblock \bibinfo{journal}{International Journal for Uncertainty
  Quantification} \bibinfo{volume}{12}.
\bibitem[{{Quintana-Segu{\'\i}} et~al.(2008){Quintana-Segu{\'\i}}, {Le Moigne},
  {Durand}, {Martin}, {Habets}, {Baillon}, {Canellas}, {Franchisteguy} and
  {Morel}}]{Quintana2008}
\bibinfo{author}{{Quintana-Segu{\'\i}}, P.}, \bibinfo{author}{{Le Moigne}, P.},
  \bibinfo{author}{{Durand}, Y.}, \bibinfo{author}{{Martin}, E.},
  \bibinfo{author}{{Habets}, F.}, \bibinfo{author}{{Baillon}, M.},
  \bibinfo{author}{{Canellas}, C.}, \bibinfo{author}{{Franchisteguy}, L.},
  \bibinfo{author}{{Morel}, S.}, \bibinfo{year}{2008}.
\newblock \bibinfo{title}{{Analysis of Near-Surface Atmospheric Variables:
  Validation of the SAFRAN Analysis over France}}.
\newblock \bibinfo{journal}{Journal of Applied Meteorology and Climatology}
  \bibinfo{volume}{47}, \bibinfo{pages}{92}.
\newblock \DOIprefix\doi{10.1175/2007JAMC1636.1}.
\bibitem[{Razavi and Gupta(2015)}]{Razavi_Gupta_SA_2014WR}
\bibinfo{author}{Razavi, S.}, \bibinfo{author}{Gupta, H.V.},
  \bibinfo{year}{2015}.
\newblock \bibinfo{title}{What do we mean by sensitivity analysis? the need for
  comprehensive characterization of “global” sensitivity in earth and
  environmental systems models}.
\newblock \bibinfo{journal}{Water Resources Research} \bibinfo{volume}{51},
  \bibinfo{pages}{3070--3092}.
\newblock \URLprefix
  \url{https://agupubs.onlinelibrary.wiley.com/doi/abs/10.1002/2014WR016527},
  \DOIprefix\doi{https://doi.org/10.1002/2014WR016527},
  \href{http://arxiv.org/abs/https://agupubs.onlinelibrary.wiley.com/doi/pdf/10.1002/2014WR016527}{{\tt
  arXiv:https://agupubs.onlinelibrary.wiley.com/doi/pdf/10.1002/2014WR016527}}.
\bibitem[{Razavi and Gupta(2019)}]{razavi2019multi}
\bibinfo{author}{Razavi, S.}, \bibinfo{author}{Gupta, H.V.},
  \bibinfo{year}{2019}.
\newblock \bibinfo{title}{A multi-method generalized global sensitivity matrix
  approach to accounting for the dynamical nature of earth and environmental
  systems models}.
\newblock \bibinfo{journal}{Environmental modelling \& software}
  \bibinfo{volume}{114}, \bibinfo{pages}{1--11}.
\bibitem[{Ross et~al.(2015)Ross, Abbey, Bouffard and
  Jos}]{ross2015multiobjective}
\bibinfo{author}{Ross, M.}, \bibinfo{author}{Abbey, C.},
  \bibinfo{author}{Bouffard, F.}, \bibinfo{author}{Jos, G.},
  \bibinfo{year}{2015}.
\newblock \bibinfo{title}{Multiobjective optimization dispatch for microgrids
  with a high penetration of renewable generation}.
\newblock \bibinfo{journal}{IEEE Transactions on Sustainable Energy}
  \bibinfo{volume}{6}, \bibinfo{pages}{1306--1314}.
\bibitem[{Roux et~al.(2011)Roux, Labat, Garambois, Maubourguet, Chorda and
  Dartus}]{roux2011physically}
\bibinfo{author}{Roux, H.}, \bibinfo{author}{Labat, D.},
  \bibinfo{author}{Garambois, P.A.}, \bibinfo{author}{Maubourguet, M.M.},
  \bibinfo{author}{Chorda, J.}, \bibinfo{author}{Dartus, D.},
  \bibinfo{year}{2011}.
\newblock \bibinfo{title}{A physically-based parsimonious hydrological model
  for flash floods in mediterranean catchments}.
\newblock \bibinfo{journal}{Natural Hazards and Earth System Sciences}
  \bibinfo{volume}{11}, \bibinfo{pages}{2567--2582}.
\bibitem[{Saharia et~al.(2021)Saharia, Kirstetter, Vergara, Gourley, Emmanuel
  and Andrieu}]{saharia2021impact}
\bibinfo{author}{Saharia, M.}, \bibinfo{author}{Kirstetter, P.E.},
  \bibinfo{author}{Vergara, H.}, \bibinfo{author}{Gourley, J.J.},
  \bibinfo{author}{Emmanuel, I.}, \bibinfo{author}{Andrieu, H.},
  \bibinfo{year}{2021}.
\newblock \bibinfo{title}{On the impact of rainfall spatial variability,
  geomorphology, and climatology on flash floods}.
\newblock \bibinfo{journal}{Water Resources Research} \bibinfo{volume}{57},
  \bibinfo{pages}{e2020WR029124}.
\bibitem[{Sahraei et~al.(2020)Sahraei, Asadzadeh and
  Unduche}]{sahraei2020signature}
\bibinfo{author}{Sahraei, S.}, \bibinfo{author}{Asadzadeh, M.},
  \bibinfo{author}{Unduche, F.}, \bibinfo{year}{2020}.
\newblock \bibinfo{title}{Signature-based multi-modelling and multi-objective
  calibration of hydrologic models: Application in flood forecasting for
  canadian prairies}.
\newblock \bibinfo{journal}{Journal of Hydrology} \bibinfo{volume}{588},
  \bibinfo{pages}{125095}.
\bibitem[{Saltelli(2002)}]{saltelli2002making}
\bibinfo{author}{Saltelli, A.}, \bibinfo{year}{2002}.
\newblock \bibinfo{title}{Making best use of model evaluations to compute
  sensitivity indices}.
\newblock \bibinfo{journal}{Computer physics communications}
  \bibinfo{volume}{145}, \bibinfo{pages}{280--297}.
\bibitem[{Saltelli et~al.(2008)Saltelli, Ratto, Andres, Campolongo, Cariboni,
  Debora~Gatelli and Tarantola}]{Saltelli_2008}
\bibinfo{author}{Saltelli, A.}, \bibinfo{author}{Ratto, M.},
  \bibinfo{author}{Andres, T.}, \bibinfo{author}{Campolongo, F.},
  \bibinfo{author}{Cariboni, J.}, \bibinfo{author}{Debora~Gatelli, M.S.},
  \bibinfo{author}{Tarantola, S.}, \bibinfo{year}{2008}.
\newblock \bibinfo{title}{Global Sensitivity Analysis: The Primer}.
\newblock \bibinfo{publisher}{Wiley, New York}.
\bibitem[{Samaniego et~al.(2010)Samaniego, Kumar and
  Attinger}]{samaniego2010multiscale}
\bibinfo{author}{Samaniego, L.}, \bibinfo{author}{Kumar, R.},
  \bibinfo{author}{Attinger, S.}, \bibinfo{year}{2010}.
\newblock \bibinfo{title}{Multiscale parameter regionalization of a grid-based
  hydrologic model at the mesoscale}.
\newblock \bibinfo{journal}{Water Resources Research} \bibinfo{volume}{46}.
\bibitem[{Shafii and Tolson(2015)}]{Shafii2014WR016520}
\bibinfo{author}{Shafii, M.}, \bibinfo{author}{Tolson, B.A.},
  \bibinfo{year}{2015}.
\newblock \bibinfo{title}{Optimizing hydrological consistency by incorporating
  hydrological signatures into model calibration objectives}.
\newblock \bibinfo{journal}{Water Resources Research} \bibinfo{volume}{51},
  \bibinfo{pages}{3796--3814}.
\newblock \URLprefix
  \url{https://agupubs.onlinelibrary.wiley.com/doi/abs/10.1002/2014WR016520},
  \DOIprefix\doi{https://doi.org/10.1002/2014WR016520},
  \href{http://arxiv.org/abs/https://agupubs.onlinelibrary.wiley.com/doi/pdf/10.1002/2014WR016520}{{\tt
  arXiv:https://agupubs.onlinelibrary.wiley.com/doi/pdf/10.1002/2014WR016520}}.
\bibitem[{Sobol and Kucherenko(2010)}]{sobol2010derivative}
\bibinfo{author}{Sobol, I.}, \bibinfo{author}{Kucherenko, S.},
  \bibinfo{year}{2010}.
\newblock \bibinfo{title}{Derivative based global sensitivity measures}.
\newblock \bibinfo{journal}{Procedia-Social and Behavioral Sciences}
  \bibinfo{volume}{2}, \bibinfo{pages}{7745--7746}.
\bibitem[{Sobol'(1990)}]{sobol1990sensitivity}
\bibinfo{author}{Sobol', I.M.}, \bibinfo{year}{1990}.
\newblock \bibinfo{title}{On sensitivity estimation for nonlinear mathematical
  models}.
\newblock \bibinfo{journal}{Matematicheskoe modelirovanie} \bibinfo{volume}{2},
  \bibinfo{pages}{112--118}.
\bibitem[{Song et~al.(2015)Song, Zhang, Zhan, Xuan, Ye and Xu}]{SONG2015}
\bibinfo{author}{Song, X.}, \bibinfo{author}{Zhang, J.}, \bibinfo{author}{Zhan,
  C.}, \bibinfo{author}{Xuan, Y.}, \bibinfo{author}{Ye, M.},
  \bibinfo{author}{Xu, C.}, \bibinfo{year}{2015}.
\newblock \bibinfo{title}{Global sensitivity analysis in hydrological modeling:
  Review of concepts, methods, theoretical framework, and applications}.
\newblock \bibinfo{journal}{Journal of Hydrology} \bibinfo{volume}{523},
  \bibinfo{pages}{739--757}.
\newblock \URLprefix
  \url{https://www.sciencedirect.com/science/article/pii/S0022169415001249},
  \DOIprefix\doi{https://doi.org/10.1016/j.jhydrol.2015.02.013}.
\bibitem[{Tarasova et~al.(2018)Tarasova, Basso, Zink and
  Merz}]{tarasova2018exploring}
\bibinfo{author}{Tarasova, L.}, \bibinfo{author}{Basso, S.},
  \bibinfo{author}{Zink, M.}, \bibinfo{author}{Merz, R.}, \bibinfo{year}{2018}.
\newblock \bibinfo{title}{Exploring controls on rainfall-runoff events: 1. time
  series-based event separation and temporal dynamics of event runoff response
  in germany}.
\newblock \bibinfo{journal}{Water Resources Research} \bibinfo{volume}{54},
  \bibinfo{pages}{7711--7732}.
\bibitem[{Tavakkoli-Moghaddam et~al.(2011)Tavakkoli-Moghaddam, Azarkish and
  Sadeghnejad-Barkousaraie}]{tavakkoli2011new}
\bibinfo{author}{Tavakkoli-Moghaddam, R.}, \bibinfo{author}{Azarkish, M.},
  \bibinfo{author}{Sadeghnejad-Barkousaraie, A.}, \bibinfo{year}{2011}.
\newblock \bibinfo{title}{A new hybrid multi-objective pareto archive pso
  algorithm for a bi-objective job shop scheduling problem}.
\newblock \bibinfo{journal}{Expert Systems with Applications}
  \bibinfo{volume}{38}, \bibinfo{pages}{10812--10821}.
\bibitem[{Torres-Trevi{\~n}o et~al.(2011)Torres-Trevi{\~n}o, Reyes-Valdes,
  L{\'o}pez and Praga-Alejo}]{torres2011multi}
\bibinfo{author}{Torres-Trevi{\~n}o, L.M.}, \bibinfo{author}{Reyes-Valdes,
  F.A.}, \bibinfo{author}{L{\'o}pez, V.}, \bibinfo{author}{Praga-Alejo, R.},
  \bibinfo{year}{2011}.
\newblock \bibinfo{title}{Multi-objective optimization of a welding process by
  the estimation of the pareto optimal set}.
\newblock \bibinfo{journal}{Expert systems with applications}
  \bibinfo{volume}{38}, \bibinfo{pages}{8045--8053}.
\bibitem[{Veluscek et~al.(2015)Veluscek, Kalganova, Broomhead and
  Grichnik}]{veluscek2015composite}
\bibinfo{author}{Veluscek, M.}, \bibinfo{author}{Kalganova, T.},
  \bibinfo{author}{Broomhead, P.}, \bibinfo{author}{Grichnik, A.},
  \bibinfo{year}{2015}.
\newblock \bibinfo{title}{Composite goal methods for transportation network
  optimization}.
\newblock \bibinfo{journal}{Expert Systems with Applications}
  \bibinfo{volume}{42}, \bibinfo{pages}{3852--3867}.
\bibitem[{Vrugt et~al.(2008)Vrugt, Ter~Braak, Clark, Hyman and
  Robinson}]{vrugt2008treatment}
\bibinfo{author}{Vrugt, J.A.}, \bibinfo{author}{Ter~Braak, C.J.},
  \bibinfo{author}{Clark, M.P.}, \bibinfo{author}{Hyman, J.M.},
  \bibinfo{author}{Robinson, B.A.}, \bibinfo{year}{2008}.
\newblock \bibinfo{title}{Treatment of input uncertainty in hydrologic
  modeling: Doing hydrology backward with markov chain monte carlo simulation}.
\newblock \bibinfo{journal}{Water Resources Research} \bibinfo{volume}{44}.
\bibitem[{Wang et~al.(2017)Wang, Zhao, Wu and Wu}]{wang2017multi}
\bibinfo{author}{Wang, N.}, \bibinfo{author}{Zhao, W.j.}, \bibinfo{author}{Wu,
  N.}, \bibinfo{author}{Wu, D.}, \bibinfo{year}{2017}.
\newblock \bibinfo{title}{Multi-objective optimization: a method for selecting
  the optimal solution from pareto non-inferior solutions}.
\newblock \bibinfo{journal}{Expert Systems with Applications}
  \bibinfo{volume}{74}, \bibinfo{pages}{96--104}.
\bibitem[{Wang and Rangaiah(2017)}]{wang2017application}
\bibinfo{author}{Wang, Z.}, \bibinfo{author}{Rangaiah, G.P.},
  \bibinfo{year}{2017}.
\newblock \bibinfo{title}{Application and analysis of methods for selecting an
  optimal solution from the pareto-optimal front obtained by multiobjective
  optimization}.
\newblock \bibinfo{journal}{Industrial \& Engineering Chemistry Research}
  \bibinfo{volume}{56}, \bibinfo{pages}{560--574}.
\bibitem[{Wu et~al.(2021)Wu, Chen, Ye, Guo, Meng and Zhang}]{wu2021improved}
\bibinfo{author}{Wu, H.}, \bibinfo{author}{Chen, B.}, \bibinfo{author}{Ye, X.},
  \bibinfo{author}{Guo, H.}, \bibinfo{author}{Meng, X.},
  \bibinfo{author}{Zhang, B.}, \bibinfo{year}{2021}.
\newblock \bibinfo{title}{An improved calibration and uncertainty analysis
  approach using a multicriteria sequential algorithm for hydrological
  modeling}.
\newblock \bibinfo{journal}{Scientific Reports} \bibinfo{volume}{11},
  \bibinfo{pages}{1--14}.
\bibitem[{Wu et~al.(2015)Wu, Liu and Lur}]{wu2015pareto}
\bibinfo{author}{Wu, Y.K.}, \bibinfo{author}{Liu, C.C.}, \bibinfo{author}{Lur,
  Y.Y.}, \bibinfo{year}{2015}.
\newblock \bibinfo{title}{Pareto-optimal solution for multiple objective linear
  programming problems with fuzzy goals}.
\newblock \bibinfo{journal}{Fuzzy Optimization and Decision Making}
  \bibinfo{volume}{14}, \bibinfo{pages}{43--55}.
\bibitem[{Yapo et~al.(1998)Yapo, Gupta and Sorooshian}]{yapo1998multi}
\bibinfo{author}{Yapo, P.O.}, \bibinfo{author}{Gupta, H.V.},
  \bibinfo{author}{Sorooshian, S.}, \bibinfo{year}{1998}.
\newblock \bibinfo{title}{Multi-objective global optimization for hydrologic
  models}.
\newblock \bibinfo{journal}{Journal of hydrology} \bibinfo{volume}{204},
  \bibinfo{pages}{83--97}.
\bibitem[{Yeh(1986)}]{yeh1986review}
\bibinfo{author}{Yeh, W.W.G.}, \bibinfo{year}{1986}.
\newblock \bibinfo{title}{Review of parameter identification procedures in
  groundwater hydrology: The inverse problem}.
\newblock \bibinfo{journal}{Water resources research} \bibinfo{volume}{22},
  \bibinfo{pages}{95--108}.
\bibitem[{Yilmaz et~al.(2008)Yilmaz, Gupta and Wagener}]{Yilmaz2007WR006716}
\bibinfo{author}{Yilmaz, K.K.}, \bibinfo{author}{Gupta, H.V.},
  \bibinfo{author}{Wagener, T.}, \bibinfo{year}{2008}.
\newblock \bibinfo{title}{A process-based diagnostic approach to model
  evaluation: Application to the nws distributed hydrologic model}.
\newblock \bibinfo{journal}{Water Resources Research} \bibinfo{volume}{44}.
\newblock \URLprefix
  \url{https://agupubs.onlinelibrary.wiley.com/doi/abs/10.1029/2007WR006716},
  \DOIprefix\doi{https://doi.org/10.1029/2007WR006716},
  \href{http://arxiv.org/abs/https://agupubs.onlinelibrary.wiley.com/doi/pdf/10.1029/2007WR006716}{{\tt
  arXiv:https://agupubs.onlinelibrary.wiley.com/doi/pdf/10.1029/2007WR006716}}.
\bibitem[{Zhu et~al.(1997)Zhu, Byrd, Lu and Nocedal}]{Zhu1997}
\bibinfo{author}{Zhu, C.}, \bibinfo{author}{Byrd, R.H.}, \bibinfo{author}{Lu,
  P.}, \bibinfo{author}{Nocedal, J.}, \bibinfo{year}{1997}.
\newblock \bibinfo{title}{Algorithm 778: L-bfgs-b: Fortran subroutines for
  large-scale bound-constrained optimization.}
\newblock \bibinfo{journal}{ACM Trans. Math. Softw.} \bibinfo{volume}{23},
  \bibinfo{pages}{550--560}.
\newblock \URLprefix
  \url{http://dblp.uni-trier.de/db/journals/toms/toms23.html#ZhuBLN97}.
\bibitem[{Zoccatelli et~al.(2011)Zoccatelli, Borga, Viglione, Chirico and
  Bl{\"o}schl}]{zoccatelli2011spatial}
\bibinfo{author}{Zoccatelli, D.}, \bibinfo{author}{Borga, M.},
  \bibinfo{author}{Viglione, A.}, \bibinfo{author}{Chirico, G.},
  \bibinfo{author}{Bl{\"o}schl, G.}, \bibinfo{year}{2011}.
\newblock \bibinfo{title}{Spatial moments of catchment rainfall: rainfall
  spatial organisation, basin morphology, and flood response}.
\newblock \bibinfo{journal}{Hydrology and Earth System Sciences}
  \bibinfo{volume}{15}, \bibinfo{pages}{3767--3783}.

\end{thebibliography}



\end{document}